\documentclass[a4paper,11pt,reqno]{amsart}
\usepackage[latin1]{inputenc}
\usepackage[english]{babel}
\usepackage{amssymb,amsfonts, amsmath, color}
\usepackage[margin=2.7cm]{geometry}
\usepackage{xcolor}
\usepackage{graphicx, color, enumerate}
\usepackage[latin1]{inputenc}
\usepackage[active]{srcltx}
\usepackage{tikz}
\usepackage{pgf}
\usepackage{etex}
\usepackage{verbatim}
\usepackage{tikz-3dplot}
\usepackage{amsmath}
\usepackage{amsfonts}
\usepackage{amssymb}
\usepackage{amsthm}
\usepackage{algpseudocode}
\usepackage{float}
\usepackage{xcolor}
\usepackage{mathrsfs}
\usepackage{import}
\usepackage{geometry}
\usepackage{fancyhdr}
\usepackage{fp}

\newtheorem{theorem}{Theorem}[section]
\newtheorem{proposition}[theorem]{Proposition}

\newtheorem{lemma}[theorem]{Lemma}
\newtheorem{remark}[theorem]{Remark}
\newtheorem{definition}[theorem]{Definition}

\newcommand{\bcl}{\begin{center}}
\newcommand{\ecl}{\end{center}}
\newcommand{\brl}{\begin{right}}
\newcommand{\erl}{\end{right}}
\newcommand{\ben}{\begin{enumerate}}
\newcommand{\een}{\end{enumerate}}
\newcommand{\overliner}{\begin{array}}
\newcommand{\earr}{\end{array}}
\newcommand{\btab}{\begin{tabular}}
\newcommand{\etab}{\end{tabular}}
\newcommand{\bdoc}{\begin{document}}
\newcommand{\edoc}{\end{document}}
\newcommand{\beqy}{\begin{eqnarray}}
\newcommand{\eeqy}{\end{eqnarray}}

\newcommand{\beqi}{\begin{eqnarray*}}
\newcommand{\eeqi}{\end{eqnarray*}}
\newcommand{\bitem}{\begin{itemize}}

\newcommand{\eitem}{\end{itemize}}
\newcommand{\nln}{\newline}
\newcommand{\newt}{\newtheorem}

\newcommand{\pa}{\partial}
\newcommand{\re}{{I\!\!R}}
\newcommand{\Rn}{\R^N}
\newcommand{\xr}{x\in\R }
\newcommand{\x}{\times}
\newcommand{\dyle}{\displaystyle}
\newcommand{\ene}{{I\!\!N}}
\newcommand{\irn}{\int\limits_{\R^N}}
\newcommand{\io}{\int\limits_{\O}}
\newcommand{\meas}{{\rm meas\,}}
\newcommand{\dif}{\nabla_{xy}}
\newcommand{\sign}{{\rm sign}}
\newcommand{\map}{\longrightarrow }
\newcommand{\imp}{\Longrightarrow }
\renewcommand{\div}{\nabla\cdot }
\newcommand{\sen}{{\rm sen\,}}
\newcommand{\tg}{{\rm tg\,}}
\newcommand{\arcsen}{{\rm arcsen\,}}
\newcommand{\arctg}{{\rm arctg\,}}
\newcommand{\supp}{{\textsl supp\ }}
\newcommand{\ity}{\int_{-\iy}^{+\iy}}
\newcommand{\limit}{\lim\limits}
\newcommand{\limi}{\limit_{n\to\infty}}
\newcommand{\sumi}{\sum\limits_{n=1}^{\infty}}
\newcommand{\ulu}{\underline u}
\newcommand{\ulw}{\underline w}
\newcommand{\ulz}{\underline z}
\newcommand{\ulv}{\underline v}
\newcommand{\uls}{\underline s}
\newcommand{\olu}{\overline u}
\newcommand{\olv}{\overline v}
\newcommand{\ols}{\overline s}
\newcommand{\ob}{\overline\b}
\newcommand{\ovar}{\overline\var}
\newcommand{\wv}{\widetilde v}
\newcommand{\wu}{\widetilde u}
\newcommand{\ws}{\widetilde s}
\renewcommand{\a }{\alpha }
\renewcommand{\b }{\beta }
\newcommand{\g }{\gamma}
\newcommand{\G }{\Gamma }
\renewcommand{\d }{\delta }

\newcommand{\D }{\Delta }
\newcommand{\e }{\varepsilon }
\newcommand{\z }{\zeta }
\renewcommand{\l }{\lambda }
\renewcommand{\L }{\Lambda }
\newcommand{\m }{\mu }
\newcommand{\n }{\nabla }
\newcommand{\s }{\sigma }
\newcommand{\Sig }{\Sigma }
\renewcommand{\t }{\tau }
\newcommand{\var }{\varphi }
\renewcommand{\o }{\omega }
\renewcommand{\O }{\Omega }
\newcommand{\R}{{\mathbb{R}}}
\newcommand{\bC}{{\bf C}}
\newcommand{\bZ}{{\bf Z}}
\newcommand{\bN}{{\bf N}}
\newcommand{\bQ}{{\bf Q}}
\newcommand{\bK}{{\bf K}}
\newcommand{\bI}{{\bf I}}
\newcommand{\bv}{{\bf v}}
\newcommand{\bV}{{\bf V}}
\DeclareMathOperator{\suppo}{supp} \DeclareMathOperator{\di}{div}

\newenvironment{Proof}{\Rmovelastskip\vskip12pt
plus 1pt \noindent\em\rm}{\hfill {\qed \hskip .2cm}}

\begin{document}

\title[]{On the uniqueness for the heat equation with density on infinite graphs}

\author{Giulia Meglioli}

\address{\hbox{\parbox{5.7in}{\medskip \noindent{Giulia Meglioli, \\Fakult\"at f\"ur Mathematik, \\Universit\"at Bielefeld, \\33501, Bielefeld, Germany \\ [3pt] \emph{E-mail address: }{\tt gmeglioli@math.uni-bielefeld.de}}}}}


\subjclass[2010]{}

\maketitle

\begin{abstract} 
We study the uniqueness of solutions to a class of heat equations with positive density posed on infinite weighted graphs. We separately consider the case when the density is bounded from below by a positive constant and the case of possibly vanishing density, showing that these two scenarios lead to two different classes of uniqueness.
\end{abstract}

\section{Introduction}\label{sec0}

In what follows we investigate uniqueness of solutions to the following initial-value problem
\begin{equation}\label{problema}
\begin{cases}
\rho\,\partial_t u-\Delta u =0 &\quad \text{in}\,\, S_T,\\
u=0&\quad \text{in}\,\,G\times\{0\},
\end{cases}
\end{equation}
where $\rho:G\to(0,+\infty)$ is a positive function usually referred to as density or weight. Furthermore we define $S_T:=G\times(0,T]$ where $(G, \omega, \mu)$ is an infinite weighted graph with edge-weight $\omega$ and {node (or vertex) measure} $\mu$ and $\Delta$ denotes the Laplace operator on $G$. We will specify later the precise assumptions on the graph $G$, see \eqref{e7f}, and on the function $\rho$, see \eqref{rhobound} and \eqref{rhobound3}. Specifically, we will demonstrate that the only function which solves \eqref{problema} is  $u\equiv0$. Consequently, we can deduce the uniqueness of solutions to the following problem:
$$
\begin{cases}
\rho\,\partial_t u-\Delta u =\,f &\quad \text{in}\,\, S_T,\\
u=u_0&\quad \text{in}\,\,G\times\{0\},
\end{cases}
$$
for some given functions $u_0$, $f$.
\medskip

Uniqueness results for the Cauchy problem given in \eqref{problema} have been thoroughly investigated in recent years, both in the context of the Euclidean setting and on Riemannian manifolds. Regarding the case when the problem \eqref{problema} is posed on $\R^N$ or on bounded domains of $\R^N$, we refer the reader to, e.g., \cite{AB,EKP,F,IKO,KPT,NP,Punzo,T}. In the case of a geodesically complete Riemannian manifold $M$, with $\Delta$ denoting the Laplace-Beltrami operator, it is worth noting that the uniqueness of bounded solutions to \eqref{problema} with $\rho \equiv 1$ is equivalent to the property of stochastic completeness of the manifold $M$; see \cite[Theorem 6.2]{Grig}. Given the heat kernel $p$ and the Riemannian volume element $d\mu$, we say that $M$ is stochastically complete if
$$
\int_M p(x,y,t)\,d\mu(y)=1\quad \text{for all}\,\,\,x\in M,\,\,t>0.
$$
We refer the reader, e.g., to \cite{Gaf,Gri,Hsu,Yau} for some results regarding the uniqueness of bounded solutions to problem \eqref{problema}. Beyond bounded solutions, the uniqueness of solutions has also been addressed in broader classes of functions, such as solutions satisfying certain integral conditions; see, e.g., \cite{Gri,Grig,L,MR,Punzo2}. In particular, in \cite{Gri,Grig}, the author proved the uniqueness of solutions to the Cauchy problem within the class of functions satisfying 
\begin{equation}\label{eq01} 
\int_0^T\int_{B_r(x_0)}u^p(x,t)\,d\mu(x)\,dt\le C\exp{f(r)} 
\end{equation} 
with $p=2$. Here, $x_0\in M$ is a fixed reference point, $B_r(x_0)$ is a geodesic ball, and $f(r)$ is a non-decreasing, positive function on $(0,+\infty)$ such that
$$
\int_1^{\infty}\frac{r}{f(r)}\,dr=\infty\,.
$$
This result was generalized in \cite{Punzo2} to the class of functions satisfying \eqref{eq01} with $p\in[1,2)$. 
\medskip

More recently, similar questions have also been addressed in the context of graphs, specifically when $G$ is an infinite weighted graph as in problem \eqref{problema}. We particularly mention reference \cite{Huang}, where problem \eqref{problema} with $\rho\equiv1$ was considered. Let $\{R_n\}$ be an increasing sequence of positive numbers with $\lim_{n\to\infty}R_n=+\infty$; the author establishes uniqueness of solutions in the class of functions satisfying, for sufficiently large $n$, the following inequality
\begin{equation}\label{eq02}
\int_0^T\sum_{x\in B_{R_n}(x_0)}u^2(x,t)\mu(x)\,dt\,\le\, C\exp\{c\,R_n\log(R_n)\},
\end{equation}
where $C>0$, $0<c<\frac12$ and $x_0\in G$ is a fixed reference point. Furthermore, the author constructs an example, in the special case of the lattice $\mathbb{Z}$, showing that uniqueness fails when solutions satisfy \eqref{eq02} with $c\ge2\sqrt2$. Another uniqueness result concerning problem \eqref{problema} (with $\rho\equiv1$) is proved in \cite{HKS}; to explain such a result, we need the following definition: we say that a graph $G$ is {\it globally local} with respect to a given monotone increasing function $f:(0,+\infty)\to(0,+\infty)$, if there exists a constant $A>1$ such that
$$
\limsup_{r\to\infty}\frac{s_r f(Ar)}{r}<\infty.
$$
Here, given a reference point $x_0\in G$, we define $s_r$ as the jump size outside $B_r$, i.e. $$s_r:=\sup\{d(x,y):\,\,x,y\in G\,\,\text{and}\,\,\,d(x,x_0)\land d(y,x_0)\ge r\} \quad\text{and}\quad s_0<\infty.$$ The authors consider problem \eqref{problema} with the additional assumption that the graph $G$ is globally local. They show that the uniqueness of bounded solutions to such a problem is equivalent to the property of stochastic completeness of $G$ .
\medskip

Similar questions have also been addressed in the elliptic counterpart of problem \eqref{problema}; see, e.g., \cite{AS2,BMP,BPgraph,MPgraph,Svagna}. More generally, PDEs posed on infinite graphs are the subject of several studies. Without aiming to be exhaustive, we mention, in the case of parabolic equations, \cite{BCG,BSW,CGZ,GMP,HM,LW2,MoPuSo2,M,PS,W} and, for elliptic equations, \cite{AS,GLY,GLY2,MoPuSo}.

\subsection{Outline of our results}
In the present paper we aim to prove that $u\equiv0$ is the only solution to problem \eqref{problema}, whenever $u$ belongs to a certain class of functions satisfying an appropriate summability condition. To do so, we distinguish between to different cases depending on the density function $\rho$. 
\smallskip

\noindent First, we assume that $\rho$ is a positive function bounded from below by a positive constant $\rho_0$ (see \eqref{rhobound}). In this case, let $p\ge 2$ and consider the class of functions $u$ such that, for two given constants $C>0$, $\beta\in\left(0,\frac 1{2s}\right)$, for some $x_0\in G$ and sufficiently large $n$, they satisfy 
\begin{equation}\label{h1}
\int_0^T\sum_{x\in B_{R_n}(x_0)}|u(x,t)|^p\,\mu(x)\,dt\,\le\, C\exp\{\beta\, R_n\log R_n\},
\end{equation}
where $\{R_n\}_{n\in\mathbb{N}}$ is an increasing sequence of positive numbers such that $\lim_{n\to\infty}R_n=+\infty $. We show that within this class, the only solution to problem \eqref{problema} is $u\equiv 0$ (see Theorem \ref{teo2}). Furthermore, we extend this result to the case of $1\le p<2$ by adding an additional assumption on the graph (see Theorem \ref{teo1}). In these two statements, we consider the same class of uniqueness as in \cite{Huang}, extending the result in two directions: we demonstrate that uniqueness persists in this class of functions when a density $\rho$ as the one given in \eqref{rhobound} is added, and when $p$ is allowed to vary over the entire range $p\ge1$, rather than considering only $\rho\equiv 1$ and $p=2$ as in \cite{Huang}. 
\smallskip

\noindent As a second investigation, we enlarge the class of density functions $\rho$ by assuming that $\rho$ may vanish as $d(x,x_0)\to\infty$, for some reference point $x_0$, according to assumption \eqref{rhobound3}, i.e. 
$$
\begin{aligned}
&\text{there exists $\rho_0>0$ such that} \\
&\rho(x) \geq  \rho_0\,[d(x,x_0) + k]^{-\sigma} \quad \text{for any} \quad x \in G\\
&\text{for some $\sigma>0$ and $k>s$}.
\end{aligned}
$$
In this case, the uniqueness class we consider is smaller than the one in the previous case. Specifically, the set of solutions satisfies the following: there exist an increasing sequence of positive numbers $\{R_n\}_{n\in\mathbb{N}}$ such that $\lim_{n\to\infty}R_n=+\infty $ and two constants $C>0$, $\beta>0$ such that, for some $x_0\in G$ and for all sufficiently large $n$
\begin{equation}\label{h2}
\int_0^T\sum_{x\in B_{R_n}(x_0)}|u(x,t)|^p\,\mu(x)\,dt\,\le\,  (R_n+k)^\beta.
\end{equation}
Then, if $p\ge2$, we show that the only solution is the trivial one if $0<\sigma\le2$ (see Theorem \ref{teo3}). Furthermore, we can extend this statement to the range of $p\in[1,2)$ if we restrict the class of graphs and if we consider the smaller interval $0<\sigma\le1$ (see Theorem \ref{teo4}).

It is worth mentioning that, even on Riemannian manifolds, results in the set $p\in[1,2)$ are more delicate and require stronger assumptions on the underlying manifold (see \cite{Punzo2}). Our results are based on suitable a priori estimates, obtained by multiplying the equation by appropriate test functions and then integrating by parts according to formula \eqref{e4f}. The difference that allows us to extend from $p\ge2$ to $p\ge1$ lies precisely in the number of times we integrate by parts. In both cases of $u$ satisfying \eqref{h1} or \eqref{h2}, performing the integration by parts just once requires $p\ge2$. On the other hand, to allow $p$ to vary in $[1,2)$ we need to integrate by part twice, which consequently necessitates stronger assumptions on the graph and, in the case of $\rho$ vanishing, also on $\sigma$.

\subsection{Organization of the paper}
The paper is organized as follows. In Section \ref{mf}, we provide an introductory overview of the graph setting along with the basic tools that will be useful throughout the paper. Section \ref{statement} is devoted to the statement of the main results. Subsequently, in Sections \ref{sec2} and \ref{funzioni}, we state and prove some auxiliary results. Finally, the main theorems are proved in Sections \ref{sec5} and \ref{sec6}, distinguishing between the case of a bounded density and a density that may vanish at infinity.

\section{Mathematical framework}\label{mf}
In this section, we gather some basic notions and definitions related to graphs, along with some preliminaries concerning analysis on graphs. All the definitions and lemmas in this section can also be found in various other works, such as \cite{BMP, GMP, Grig2, Huang, HKS, MPgraph, Svagna}.

\subsection{The graph setting}
Let $G$ be a given countable set of {\it vertices} or {\it nodes}. We introduce on $G$ the {\it edge weight} function $\omega:G\times G\to [0, +\infty)$ which satisfies the following properties:
\begin{align*}
&\text{1.} \,\,\,\omega(x,x)=0 \quad \quad \quad\,\,\text{for all $x\in G$, ({\it no loops})};\\
&\text{2.}\,\,\,\omega(x,y)=\omega(y,x) \quad\text{for all $(x,y)\in G\times G$, ({\it symmetry})};\\
&\text{3.}\,\, \displaystyle \sum_{y\in G} \omega(x,y)<\infty \quad\,\, \text{for all $x\in G$,  ({\it finite sum})}.
\end{align*}
Moreover, let $\mu:G\to (0, +\infty)$ be a given function, which we refer to as the {\it node (or vertex) measure}. 
\begin{equation*}\text{Then, a graph is the triplet $(G, \omega, \mu)$.}\end{equation*}

\noindent Additionally, we say that
\begin{itemize}
\item  The graph is {\it finite} if $G$ is a finite set; conversely, if $G$ is infinite and countable, we say that the graph is {\it infinite};
\item The pair $(x,y)$ is an {\it edge} of the graph, and we write $x\sim y$ whenever $\omega(x,y)>0$. In this case, we say that $x$ is {\it connected} to $y$ (or $x$ is {\it joined} to $y$, or $x$ is {\it adjacent} to $y$, or $x$ is a {\it neighbor} of $y$). The vertices $x$ and $y$ are called the {\it endpoints} of the edge $(x,y)$;
\item The graph is {\it undirected} if the edges do not have an orientation;
\item The graph does not posses {\it loops} if no edge connects a vertex to itself, and no two vertices are connected by more than one edge; 
\item A collection of vertices $\gamma\equiv \{x_k\}_{k=0}^n\subset G$ is a {\it path} if $x_k\sim x_{k+1}$ for all $k=0, \ldots, n-1.$ A graph is {\it connected} if, for any two vertices $x,y\in V$, there exists a path joining $x$ to $y$. 
\item  A graph  is {\it locally finite}, if, for each vertex $x\in G$, $\operatorname{card}\{y\in G\,:\, x\sim y\}<\infty\,.$
\end{itemize}
We define the {\it degree} $\operatorname{deg}$ and the {\it weighted degree} $\operatorname{Deg}$ of a vertex $x\in G$ as follows:
\[\operatorname{deg}(x):=\sum_{y\in G}\omega(x,y), \quad \operatorname{Deg}(x):=\frac{\operatorname{deg}(x)}{\mu(x)}\,.\]

\noindent A {\it pseudo metric} on $G$ is a mapping $d:G\times G\to [0, +\infty)$ such that
\begin{itemize}
  \item $d(x,x)=0$ \quad \text{ for all } $x\in G$;
  \item $d(x,y)=d(y,x)$ \text{ for all } $x,y\in G$;
  \item $d(x,y)\leq d(x,z)+d(z,y)$ \text{ for all } $x,y,z\in G$.
\end{itemize}
In general, $d$ is not a metric, since we can find points $x, y\in G$, $x\neq y$ such that $d(x,y)=0\,.$ We conclude this subsection with two fundamental definitions.

\begin{definition}
We define the {\em jump size} $s$ of the pseudo metric $d$ as
\begin{equation}\label{e14f}
s:=\sup\{d(x,y) \,:\, x,y\in G, \omega(x,y)>0\}.
\end{equation}
\end{definition}

\begin{definition}\label{intrmet}
Let $C_0>0.$ We say that a pseudo metric $d$ on $(G, \omega, \mu)$ is {\em$q-$intrinsic} with bound $C_0$ if
\begin{equation}\label{e1f}
\frac 1{\mu(x)}\sum_{y\in G} \omega(x,y)d^q(x,y) \leq C_0\quad \text{ for all } x\in G\,.
\end{equation}
If $q=2, C_0=1$, then $d$ is usually called {\em intrinsic} (see e.g. \cite{HKS}).
\end{definition}

Hypothesis \eqref{e1f} has an analogous counterpart on Riemannian manifolds. We refer the reader to \cite[Remark 2.4]{MPgraph} for a specific comparison between the two settings. 

\subsection{Difference and Laplace operators on graphs} Let $\mathfrak F$ denote the set of all functions $f: G\to \mathbb R$\, and let $ \mathfrak F_T$ denote the set of all functions $f:G\times[0,T]\to \R$. For any $f\in \mathfrak F$, we introduce the {\it difference operator}
\begin{equation}\label{e2f}
\nabla_{xy} f:= f(y)-f(x)\,.
\end{equation}
It is direct to see that for any $f,g\in \mathfrak F$, it is possible to derive the following {\it product rule} 
\begin{equation}\label{e3f}
\nabla_{xy}(fg)=f(x) (\nabla_{xy} g) + (\nabla_{xy} f)g(y) \quad \text{ for all } x,y\in G\,.
\end{equation}
The gradient squared of $f\in \mathfrak F$ is defined as (see \cite{CGZ})
\[|\nabla f(x)|^2=\frac 1{\mu(x)}\sum_{y\in G}\omega(x,y)(\dif f)^2, \quad x\in G\,.\]

\begin{definition}\label{def1}
Let $(G, \omega,\mu)$ be a weighted graph. For any $f\in \mathfrak F$ such that 
\begin{equation}\label{sumfinite}
\sum_{y\in G}\omega(x, y) |f(y)|<\infty \quad \text{ for all } \; x\in G,
\end{equation}
 the {\em (weighted) Laplace operator} of $(G, \omega, \mu)$ is
\begin{equation}\label{eq11}
\Delta f(x):=\frac{1}{\mu(x)}\sum_{y\in G}[f(y)-f(x)]\omega(x,y)\quad \text{ for all }\, x\in G\,.
\end{equation}
\end{definition}
\noindent Note that \eqref{sumfinite} holds whenever the graph is locally finite. Clearly,
\[\Delta f(x)=\frac 1{\mu(x)}\sum_{y\in G}(\dif f)\omega(x,y)\quad \text{ for all } x\in G\,.\]

It is easily seen that if at least one of the functions $f, g\in \mathfrak F$ has {\it finite} support, then the following {\it integration by parts formula} holds:
\begin{equation}\label{e4f}
\sum_{x\in G}[\Delta f(x)] g(x) \mu(x)=-\frac 1 2\sum_{x,y\in G}(\dif f)(\dif g) \omega(x,y)\,.
\end{equation}

\noindent We also recall the following two lemmas which will be expedient in the sequel. For the proofs of such results we refer the reader to \cite{MPgraph}.

\begin{lemma}\label{lemma4} 
Let $\psi\in C^1(\mathbb R; \mathbb R)$ be a convex function, let $u\in \mathfrak F$. Then
\begin{equation}\label{e5f}
\Delta \psi(u(x)) \geq \psi'(u(x)) \Delta u(x)\quad \text{ for all }\; x\in G\,.
\end{equation}
\end{lemma}

\begin{lemma}[\textit{Laplacian of the product}]\label{lemma5} 
Let $f,g\in \mathfrak F$. Then, for every $x\in G$,
\begin{equation}\label{e25f}
\Delta [f(x)g(x)]=f(x) \Delta g(x) + g(x) \Delta f(x) + \frac 1{\mu(x)}\sum_{y\in G}(\dif f)(\dif g )\omega(x,y)\,.
\end{equation}
\end{lemma}

\section{Statement of the main results}\label{statement}\setcounter{equation}{0}

For any $x_0\in G$ and $r>0$ we set
\[B_r(x_0):=\{x\in G\,:\, d(x,x_0)<r\}\,.\] 
In what follows, we will always make the following assumption:
\begin{equation}\label{e7f}
\begin{cases}
(i) & (G, \mu, \omega) \text{ is an infinite, connected, weighted graph};\\
(ii) & \text{there exists an intrinsic pseudo metric } d \text{ such that  } \\
& \text{the jump size $s$ is finite}\\
(iii)& \text{the set } B_r(x_0) \text{ is a finite set,  for any } x_0\in G, r>0.\\
\end{cases}
\end{equation}
Note that by hypothesis \eqref{e7f}, $G$ is automatically locally finite. 
\medskip

\noindent Here and below, when dealing with a solution to \eqref{problema} we mean the following:
for a given $T>0$, we say that $u:S_T\to\R$ is a solution to \eqref{problema} if
\begin{itemize}
\item for every $t\in(0,T]$,
$$
\sum_{y\in G} |u(y,t)|\omega(x, y)<\infty \quad \text{ for all } \; x\in G;
$$
\item for every $x\in G$, the function $t\mapsto u(x,t)$ is continuously differentiable and $u(x,t)\to0$ as $t\to0^+$;
\item $u$ solves
$$
\rho\,\partial_t u-\Delta u =0 \quad \text{on}\,\, S_T.
$$
\end{itemize}
\medskip

\noindent We distinguish the statements of the main theorems between two different cases. First, in Subsection \ref{bounded}, we consider a density function $\rho:G\to(0,+\infty)$ that is bounded from below by a positive constant. Specifically, we assume that
\begin{equation}\label{rhobound}
\begin{aligned}
&\text{there exists $\rho_0>0$ such that} \\
&\rho(x)\ge\rho_0 \quad\text{for all}\,\,\,x\in G.
\end{aligned}
\end{equation}
Afterward, in Subsection \ref{vanishing}, we allow the positive density $\rho$ to vanish at infinity. More precisely, let $x_0\in G$ be fixed; we will assume that
\begin{equation}\label{rhobound3}
\begin{aligned}
&\text{there exists $\rho_0>0$ such that} \\
&\rho(x) \geq  \rho_0\,[d(x,x_0) + k]^{-\sigma} \quad \text{for any} \,\, x \in G\\
&\text{for some $\sigma>0$ and $k>\max\{1,s\}$},
\end{aligned}
\end{equation}

\noindent where $s$ is the jump-size defined in \eqref{e14f}.

\subsection{Density bounded away from zero}\label{bounded}
As already mentioned, the first two results deal with a weight function $\rho:G\to(0,+\infty)$ which satisfies \eqref{rhobound}.

\begin{theorem}\label{teo2}
Let $(G, \mu, \omega)$ be a weighted graph as described in \eqref{e7f}. Suppose $u$ is a solution to problem \eqref{problema}, let $p\ge2$ and assume that \eqref{rhobound} holds.
Furthermore, suppose that there exist an increasing sequence of positive numbers $\{R_n\}_{n\in\mathbb{N}}$ such that $$\lim_{n\to\infty}R_n=+\infty, $$ and two constants $C>0$, $\beta\in\left(0,\frac 1{2s}\right)$ such that, for some $x_0\in G$ and for all $n$ large enough
\begin{equation}\label{maincondition2}
\int_0^T\sum_{x\in B_{R_n}(x_0)}|u(x,t)|^p\,\mu(x)\,dt\,\le\, C\exp\{\beta\, R_n\log R_n\}.
\end{equation}
Then $u(x,t)=0$ on $S_T$.
\end{theorem}


\begin{theorem}\label{teo1}
Let $(G, \mu, \omega)$ be a weighted graph as described in \eqref{e7f}, and let the pseudo metric $d$ be $1-$intrinsic according to Definition \ref{intrmet}. Suppose $u$ is a solution to problem \eqref{problema}, let $p\geq 1$ and assume that \eqref{rhobound} holds. Furthermore, suppose that there exist an increasing sequence of positive numbers $\{R_n\}_{n\in\mathbb{N}}$ such that
$$
\lim_{n\to\infty}R_n=+\infty
$$
and two constants $C>0$, $\beta\in\left(0,\frac 1{s}\right)$ such that, for some $x_0\in G$ and for all $n$ large enough, \eqref{maincondition2} is satisfied.
Then $u(x,t)=0$ on $S_T$.
\end{theorem}

\begin{remark}
Observe that the Theorem \ref{teo2} allows us to deal only with $p\ge2$. In order to let $p$ be also smaller than $2$ we need to furthermore assume that the pseudo metric $d$ is $1-$intrinsic, see Definition \ref{intrmet}. Such a case is considered in Theorem \ref{teo1} where indeed $p\ge1$. For this reason, both theorems are interesting to be stated.
\end{remark}

\subsection{Vanishing density}\label{vanishing}
So far, we have assumed that the density function $\rho$ is bounded below by a positive constant, $\rho_0$. However, it is also interesting to consider a more general class of density functions, i.e. $\rho:G\to(0,+\infty)$ which may vanish when the distance function from a given reference point tends to infinity. In this subsection we will assume that $\rho$ satisfies \eqref{rhobound3}. Expanding the class of densities, affect the class of uniqueness we can consider, as discussed in Theorems \ref{teo3} and \ref{teo4}. For a deeper understanding of the differences, we refer the reader to Remarks \ref{remcomparison} and \ref{remcomparison2}.

\begin{theorem}\label{teo3}
Let $(G, \mu, \omega)$ be a weighted graph as described in \eqref{e7f}. Suppose $u$ is a solution to problem \eqref{problema}. Furthermore, let $p\geq 2$, assume that \eqref{rhobound3} holds with $0<\sigma \le2$ and suppose that there exist an increasing sequence of positive numbers $\{R_n\}_{n\in\mathbb{N}}$ such that $$\lim_{n\to\infty}R_n=+\infty $$ and two constants $C>0$, $\beta>0$ such that for all $n$ large enough
\begin{equation}\label{maincondition3}
\int_0^T\sum_{x\in B_{R_n}(x_0)}|u(x,t)|^p\,\mu(x)\,dt\,\le\, C\, (R_n+k)^\beta.
\end{equation}
Then $u(x,t)=0$ on $S_T$.
\end{theorem}

\begin{remark}\label{exponent}
We observe that the exponent $\sigma=2$, which arises as an upper-bound in the assumption on $\rho$, see \eqref{rhobound3} and Theorem \ref{teo3}, might be considered not a surprise. For similar results see \cite{AB,BPgraph,EKP,IKO,KPT,PT,Svagna}.
Optimality in the framework of parabolic problem as \eqref{problema} will be addressed in a forthcoming paper.
\end{remark}

\begin{theorem}\label{teo4}
Let $(G, \mu, \omega)$ be a weighted graph as described in \eqref{e7f}, and let the pseudo metric $d$ be $1-$intrinsic according to Definition \ref{intrmet}. Suppose $u$ is a solution to problem \eqref{problema}, let $p\geq 1$, assume that \eqref{rhobound3} holds with $0<\sigma \le1$ and suppose there exist an increasing sequence of positive numbers $\{R_n\}_{n\in\mathbb{N}}$ such that $$\lim_{n\to\infty}R_n=+\infty $$ and two constants $C>0$, $\beta>0$ such that for all $n$ large enough, \eqref{maincondition3} is satisfied.
Then $u(x,t)=0$ on $S_T$.
\end{theorem}

\begin{remark}
Similarly to Subsection \ref{bounded}, also in this case, the first result (i.e. Theorem \ref{teo3}) follows only for $p\ge2$. In order to cover all the set of $p\ge1$, one needs to require a $1-$intrinsic pseudo metric $d$ (see Theorem \ref{teo4}). Moreover, in this case, another restriction is necessary. Let us recall the assumption on $\rho$ introduces in \eqref{rhobound3}, i.e. $$\rho(x) \geq  \rho_0\,[d(x,x_0) + k]^{-\sigma},$$ for any $x\in G$ and for some $x_0\in G$, $k>\max\{1,s\}$, $\rho_0>0$, $\sigma>0$. In Theorem \ref{teo3}, we may let $\sigma$ lie in the interval $(0,2]$ while in Theorem \ref{teo4}, $\sigma$ must be contained in the smaller interval $(0,1]$.
\end{remark}

\begin{remark}\label{ellp}
Theorems \ref{teo2} and \ref{teo1} take into account that the density function $\rho$ is bounded from below by the positive constant $\rho_0$, as better specified in \eqref{rhobound}. Here the class of uniqueness is given by the set of functions which satisfies \eqref{maincondition2}. To better understand the situation, we observe the following. If $u\in L^p((0,T),\ell^p_{\varphi_\beta}(G, \mu))$ with 
$$
\varphi_\beta(x):=\exp\{-\beta\, d(x_, x_0)\log (d(x_, x_0))\},\quad\beta\in\left(0,\frac1{2s}\right),
$$
then, for some $C>0$, \eqref{maincondition2} follows. In fact, $u\in  L^p((0,T),\ell^p_{\varphi_\beta}(G, \mu))$ if and only if, for some $C>0$,
$$\int_0^T \sum_{x\in G}|u(x,t)|^p \mu(x) e^{-\beta d(x_, x_0)\log (d(x_, x_0))}\,dt\leq C.$$ 
Hence, for $n$ sufficiently large so that $R_n>1$, we get
\[C\geq \int_0^T\sum_{x\in B_{R_n}(x_0)}|u(x,t)|^p \mu(x) e^{-\beta d(x_, x_0)\log (d(x_, x_0))}\,dt\geq e^{-\beta R_n\log R_n}\int_0^T\sum_{x\in B_{R_n}(x_0)}|u(x,t)|^p \mu(x) \,dt\,. \]
Therefore, \eqref{maincondition2} follows. Consequently, from Theorem \ref{teo2}, it follows that the unique solution in the space $ L^p((0,T),\ell^p_{\varphi_\beta}(G, \mu))$ with $p\ge2$ is $u\equiv0$. The same result is achieved for all $p\ge1$ if we restrict the environment to those graph with a $1-$intrinsic pseudo metric and we apply Theorem \ref{teo1}.
\end{remark}

\begin{remark}\label{ellp2}
We can observe the same consequence of Remark \ref{ellp} in the case of Theorems \ref{teo3} and \ref{teo4}. Here, we assume that the density function $\rho$ vanishes as the distance from a reference point $x_0\in G$ diverges (see \eqref{rhobound3}). The class of uniqueness is given by the set of functions which satisfies \eqref{maincondition3}. Similarly to the previous case, we observe that, if $u\in L^p((0,T),\ell^p_{\psi_\beta}(G, \mu))$ with 
$$
\psi_\beta(x):=(d(x_, x_0)+k)^{-\beta},\quad\beta>0,\quad k>\max\{1,s\},
$$
then, for some $C>0$, \eqref{maincondition3} follows. Therefore, from Theorem \ref{teo3}, it follows that the unique solution to problem \eqref{problema} with $0<\sigma\le2$ in \eqref{rhobound3}, in the space $ L^p((0,T),\ell^p_{\psi_\beta}(G, \mu))$ with $p\ge2$ is $u\equiv0$. The same result is achieved, following Theorem \ref{teo1}, for all $p\ge1$ if we restrict the environment to those graph with a $1-$intrinsic pseudo metric and to the range $0<\sigma\le1$ in \eqref{rhobound3} as showed in Theorem \ref{teo4}.
\end{remark}

We conclude this section with a comparison between the results presented in Subsection \ref{bounded} and those in Subsection \ref{vanishing}. In general, the uniqueness classes of Subsections \ref{bounded} and \ref{vanishing} include, in both cases, {\it bounded} solutions. Thus, in particular, we get uniqueness of bounded solutions. 
Furthermore, also unbounded solutions can be considered but the set of unbounded functions that satisfy \eqref{maincondition2} is larger than the set of functions satisfying \eqref{maincondition3}. Moreover, the class of graphs that we can deal with in the case of Theorems \ref{teo2} and \ref{teo1} is larger than the one of Theorems \ref{teo3} and \ref{teo4}. We compare this two sets giving two examples in the forthcoming remarks.

\begin{remark}\label{remcomparison}
Let us assume that
$$
\mu(B_{R}(x_0)):=\sum_{x\in B_{R}(x_0)}\mu(x)\,\le\,CR^q\quad\text{for all}\,\,R\ge2\,\,\text{and for some}\,\,q>0.
$$
Then condition \eqref{maincondition2} is granted if, for some $g\in L^p(0,T)$ one has
$$
u(x,t)\le g(t)\exp\{\gamma d(x,x_0)\log(d(x,x_0))\}\quad\text{for all}\,\,\,(x,t)\in S_T,
$$
provided that $\gamma<\frac{\beta}{p}$. On the other hand, condition \eqref{maincondition3} follows if, for some $g\in L^p(0,T)$, one has
$$
u(x,t)\le g(t)\,(d(x,x_0)+1)^{\tau}\quad\text{for all}\,\,\,(x,t)\in S_T,,
$$
provided that $\tau<\frac{\beta-q}{p}$ and $\beta>q$.
\end{remark}

\begin{remark}\label{remcomparison2}
Let us assume that
$$
\mu(B_{R}(x_0)):=\sum_{x\in B_{R}(x_0)}\mu(x)\,\le\,Ce^{qR}\quad\text{for all}\,\,R\ge2\,\,\text{and for some}\,\,q>0.
$$
Then condition \eqref{maincondition2} is granted if, for some $g\in L^p(0,T)$ one has
$$
u(x,t)\le g(t)\exp\{\gamma d(x,x_0)\log(d(x,x_0))\}\quad\text{for all}\,\,\,(x,t)\in S_T,
$$
provided that $\gamma<\frac{\beta-q}{p}$ and $\beta>q$. Furthermore, we observe that Theorems \ref{teo3} and \ref{teo4} do not apply to this set of graphs.
\end{remark}

%

\section{Auxiliary results on infinite graphs} \label{sec2}\setcounter{equation}{0}

Before delving into the details of the proofs of the main results, we need to establish some preliminary facts that will play a key role in the following subsections. We note that Proposition \ref{prop41} is essential in the proofs of both Theorems \ref{teo2} and \ref{teo3}. On the other hand, Theorems \ref{teo1} and \ref{teo4} share the auxiliary result stated in Proposition \ref{ineqv}.

\noindent This distinction arises because Theorems \ref{teo2} and \ref{teo3} are valid for any $p\ge2$, so the result in both cases is obtained by performing only one integration by parts. Conversely, Theorems \ref{teo1} and \ref{teo4} are proved by performing a double integration by parts, allowing us to extend $p$ to the range $[1,+\infty)$.

\begin{proposition}\label{prop41}
Let $u$ be a solution to equation \eqref{problema}. Let $\eta\in \mathfrak F$, $\zeta\in\mathfrak F_T$ such that
\begin{align}
&\text{$\bullet$}\quad \eta\ge0, \; \operatorname{supp} \eta \text{ is finite };\nonumber
\\
&\text{$\bullet$}\quad \zeta\ge0,\,\,\zeta \text{ is continuously differentiable in } t \text{ on } [0,T], \text { for each } x\in G; \nonumber\\
&\text{$\bullet$}\quad [\eta^2(y)-\eta^2(x)][\zeta(y,t)-\zeta(x,t)]\ge0 \quad \text{ for all } x,y\in G, x\sim y \text{ and } t\in[0,T].\label{eq41}
\end{align}
Then, for any $p\ge 2$ and any $\tau\in (0,T]$,
\begin{equation*}\label{eq42}
\begin{aligned}
\sum_{x\in G}&|u(x,\tau)|^p\rho(x)\eta^2(x)\zeta(x,\tau)\mu(x)\le2\int_0^\tau \sum_{x,y\in G}|u(x,t)|^p\zeta(y,t)\left[\eta(y)-\eta(x)\right]^2\omega(x,y)\,dt\\
&+\int_0^\tau \sum_{x\in G}|u(x,t)|^p\eta^2(x)\left\{\rho(x)\partial_t\zeta(x,t)\mu(x)+\frac 12\zeta(x,t)\sum_{y\in G}\omega(x,y)\left[1-\frac{\zeta(y,t)}{\zeta(x,t)}\right]^2 \right\}\,dt
\end{aligned}
\end{equation*}
\end{proposition}

\begin{proof}
For any $\alpha>0$, we define
\begin{equation}\label{eq43}
\pi_{\alpha}(u):=(u^2+\alpha)^{\frac p4}\,.
\end{equation}
Now, for any $x\in G$, we consider the expression
\begin{equation}\label{eq44}
\rho(x)\,\partial_t\pi_{\alpha}(u(x,t))-\Delta \pi_{\alpha}(u(x,t)).
\end{equation}
We multiply \eqref{eq44} by $\pi_{\alpha}(u(x,t))\,\eta^2(x)\zeta(x,t)\mu(x)$ and then we sum over $x\in G$. Thus
\begin{equation}\label{eq45}
\begin{aligned}
\sum_{x\in G}\,\rho(x)\,&\partial_t \pi_{\alpha}(u(x,t))\pi_{\alpha}(u(x,t))\,\eta^2(x){\zeta(x,t)}\mu(x)\,\\
&-\sum_{x\in G}\Delta \pi_{\alpha}(u(x,t))\pi_{\alpha}(u(x,t))\,\eta^2(x){\zeta(x,t)}\mu(x)\,.\end{aligned}
\end{equation}
Note that since $\eta(x)$ is finitely supported, the series in \eqref{eq45} are finite sums. Set
\begin{equation}\label{eq46}
I:=\sum_{x\in V}\Delta \pi_{\alpha}(u(x,t))\pi_{\alpha}(u(x,t))\,\eta^2(x){\zeta(x,t)}\mu(x)\,.
\end{equation}
Then, by \eqref{e3f} and \eqref{e4f},
\begin{equation}\label{eq47}
\begin{aligned}
I&=-\frac 12 \sum_{x,y\in G}\left(\nabla_{xy} \pi_{\alpha}(u)\right) \nabla_{xy}\left[\pi_{\alpha}(u)\eta^2{\zeta}\right]\omega(x,y)\\
&=-\frac 12 \sum_{x,y\in G}\left[\nabla_{xy}\pi_{\alpha}(u)\right]^2\eta^2(y){\zeta(y,t)}\,\omega(x,y)\\
&\,\,\,\,\,-\frac 12 \sum_{x,y\in G}\pi_{\alpha}(u(x,t))\eta^2(x)\left(\nabla_{xy}\pi_{\alpha}(u)\right)\left(\nabla_{xy}{\zeta}\right)\,\omega(x,y)\\
&\,\,\,\,\,-\frac 12 \sum_{x,y\in G}\pi_{\alpha}(u(x,t)){\zeta(y,t)}\left(\nabla_{xy}\pi_{\alpha}(u)\right)\left(\nabla_{xy}\eta^2\right)\,\omega(x,y)\\
&=:J_1+J_2+J_3\,.
\end{aligned}
\end{equation}
In view of \eqref{e2f}, we obviously have
\begin{equation}\label{eq48}
\dif \eta^2=[\eta(y)+\eta(x)][\eta(y)-\eta(x)] \quad \text{ for all }\, x,y\in G\,.
\end{equation}
Due to \eqref{eq48}, by Young's inequality with exponent $2$, we have, for any $\delta_1>0$,
\begin{equation}\label{eq49}
\begin{aligned}
J_3&= -\frac 12 \sum_{x,y\in G}\pi_{\alpha}(u(x,t)){\zeta(y,t)}(\nabla_{xy}\pi_{\alpha}(u))\left[\eta(y)+\eta(x)\right]\left[\eta(y)-\eta(x)\right]\,\omega(x,y)\\
&\le \frac{\delta_1}{4} \sum_{x,y\in G}{\zeta(y,t)}\left[\nabla_{xy}\pi_{\alpha}(u)\right]^2\left[\eta(y)+\eta(x)\right]^2\,\omega(x,y)\\
&\,\,\,\,\,+\frac 1{4\delta_1} \sum_{x,y\in G}{\zeta(y,t)} \pi_{\alpha}^2(u(x,t))\left[\eta(y)-\eta(x)\right]^2\,\omega(x,y)\\
&\le \frac{\delta_1}{2} \sum_{x,y\in G}{\zeta(y,t)}\left[\nabla_{xy}\pi_{\alpha}(u)\right]^2 \left[\eta^2(y)+\eta^2(x)\right]\,\omega(x,y)\\
&\,\,\,\,\,+\frac 1{4\delta_1} \sum_{x,y\in G}{\zeta(y,t)} \pi_{\alpha}^2(u(x,t))\left[\eta(y)-\eta(x)\right]^2\,\omega(x,y)\\
&= \frac{\delta_1}{4} \sum_{x,y\in G}\left[{\zeta(y,t)}+{\zeta(x,t)}\right]\left[\nabla_{xy}\pi_{\alpha}(u)\right]^2\left[\eta^2(y)+\eta^2(x)\right]\,\omega(x,y)\\
&\,\,\,\,\,+\frac 1{4\delta_1} \sum_{x,y\in G}{\zeta(y,t)} \pi_{\alpha}^2(u(x,t))\left[\eta(y)-\eta(x)\right]^2\,\omega(x,y)\,.\\
\end{aligned}
\end{equation}
From \eqref{eq41} we can easily infer that
\begin{equation}\label{eq410}
\left[\eta^2(y)+\eta^2(x)\right]\left[{\zeta(y,t)}+{\zeta(x,t)}\right]\le 2\left\{\eta^2(y){\zeta(y,t)}+\eta^2(x){\zeta(x,t)}\right\}\,.
\end{equation}
By using \eqref{eq49} and \eqref{eq410},
\begin{equation}\label{eq411}
\begin{aligned}
J_3&\le \frac{\delta_1}{2} \sum_{x,y\in G}\left[\nabla_{xy}\pi_{\alpha}(u)\right]^2\left\{\eta^2(y){\zeta(y,t)}+\eta^2(x){\zeta(x,t)}\right\}\,\omega(x,y)\\
&\,\,\,\,\,+\frac 1{4\delta_1} \sum_{x,y\in G}{\zeta(y,t)} \pi_{\alpha}^2(u(x,t))\left[\eta(y)-\eta(x)\right]^2\,\omega(x,y)\\
&=\delta_1\sum_{x,y\in G}\left[\nabla_{xy}\pi_{\alpha}(u)\right]^2\eta^2(y){\zeta(y,t)}\,\omega(x,y)\\
&\,\,\,\,\,+\frac 1{4\delta_1} \sum_{x,y\in G}{\zeta(y,t)} \pi_{\alpha}^2(u(x,t))\left[\eta(y)-\eta(x)\right]^2\,\omega(x,y)\,.\\
\end{aligned}
\end{equation}
Similarly, by Young's inequality with exponent $2$, for any $\delta_2>0$, we have
\begin{equation}\label{eq412}
\begin{aligned}
J_2&=-\frac 12 \sum_{x,y\in G}\pi_{\alpha}(u(x,t))\eta^2(x)(\nabla_{xy}\pi_{\alpha}(u))(\nabla_{xy}{\zeta})\,\omega(x,y)\\
&=-\frac 12 \sum_{x,y\in G}\pi_{\alpha}(u(y,t))\eta^2(y)(\nabla_{xy}\pi_{\alpha}(u))\left[{\zeta(y,t)}-{\zeta(x,t)}\right]\,\omega(x,y)\\
&=-\frac 12 \sum_{x,y\in G}\pi_{\alpha}(u(y,t))\eta^2(y)(\nabla_{xy}\pi_{\alpha}(u)){\zeta(y,t)}\left[1-\frac{\zeta(x,t)}{\zeta(y,t)}\right]\,\omega(x,y)\\
&\le \frac{\delta_2}{4}\sum_{x,y\in G}\left[(\nabla_{xy}\pi_{\alpha}(u))\right]^2\eta^2(y){\zeta(y,t)}\,\omega(x,y)\\
&\quad + \frac1{4\delta_2}\sum_{x,y\in G}\pi^2_{\alpha}(u(y,t))\left[1-\frac{\zeta(x,t)}{\zeta(y,t)}\right]^2\eta^2(y){\zeta(y,t)}\,\omega(x,y)\,.
\end{aligned}
\end{equation}
By combining together \eqref{eq46}, \eqref{eq411} and \eqref{eq412} we get
\begin{equation}\label{eq413}
\begin{aligned}
\sum_{x\in G}&\Delta \pi_{\alpha}(u(x,t))\pi_{\alpha}(u(x,t))\,\eta^2(x){\zeta(x,t)}\mu(x)\\
&\le -\frac 12 \sum_{x,y\in G}[\nabla_{xy}\pi_{\alpha}(u)]^2\eta^2(y){\zeta(y,t)}\omega(x,y) \\
&\quad + \frac{\delta_2}{4}\sum_{x,y\in G}[\nabla_{xy}\pi_{\alpha}(u)]^2\eta^2(y){\zeta(y,t)}\omega(x,y)\\
&\quad + \frac1{4\delta_2}\sum_{x,y\in G}\pi^2_{\alpha}(u(y,t))\left[1-\frac{\zeta(x,t)}{\zeta(y,t)}\right]^2\eta^2(y){\zeta(y,t)}\,\omega(x,y) \\
&\quad + \delta_1\sum_{x,y\in G}\left[\nabla_{xy}\pi_{\alpha}(u)\right]^2\eta^2(y){\zeta(y,t)}\,\omega(x,y)\\
&\quad +\frac 1{4\delta_1} \sum_{x,y\in G} \pi_{\alpha}^2(u(x,t))\left[\eta(y)-\eta(x)\right]^2{\zeta(y,t)}\,\omega(x,y)\,.
\end{aligned}
\end{equation}
We now choose $\delta_1=\frac 14$ and $\delta_2=1$ in such a way that $-\frac 12+ \frac{\delta_2}{4}+ \delta_1=0$. Consequently, \eqref{eq413} yields
\begin{equation}\label{eq414}
\begin{aligned}
\sum_{x\in G}\Delta \pi_{\alpha}&(u(x,t))\pi_{\alpha}(u(x,t))\,\eta^2(x){\zeta(x,t)}\mu(x)\\
&\le \frac1{4}\sum_{x,y\in G}\pi^2_{\alpha}(u(x,t))\left[1-\frac{\zeta(y,t)}{\zeta(x,t)}\right]^2 \eta^2(x){\zeta(x,t)}\,\omega(x,y) \\
&\quad + \sum_{x,y\in G}\pi_{\alpha}^2(u(x,t))\left[\eta(y)-\eta(x)\right]^2{\zeta(y,t)}\,\omega(x,y)\,.
\end{aligned}
\end{equation}
We know treat the first term of \eqref{eq45}. We get
\begin{equation}\label{eq415}
\begin{aligned}
\sum_{x\in G}\,\rho(x)\,&\partial_t \pi_{\alpha}(u(x,t))\pi_{\alpha}(u(x,t))\,\eta^2(x){\zeta(x,t)}\mu(x)\\
&=\frac12\sum_{x\in G}\,\rho(x)\,\partial_t \pi_{\alpha}^2(u(x,t))\,\eta^2(x){\zeta(x,t)}\mu(x)\\
&=\frac12\sum_{x\in G}\,\rho(x)\,\partial_t \left[\pi_{\alpha}^2(u(x,t)){\zeta(x,t)}\right]\eta^2(x)\mu(x)\\
&\quad -\frac12\sum_{x\in G}\,\rho(x)\,\partial_t \zeta(x,t)\pi_{\alpha}^2(u(x,t))\,\eta^2(x)\mu(x)\\
\end{aligned}
\end{equation}
By substituting \eqref{eq414} and \eqref{eq415} into \eqref{eq45}, we get
\begin{equation}\label{eq416}
\begin{aligned}
\sum_{x\in G}\,&\left[\rho(x)\,\partial_t \pi_{\alpha}(u(x,t))-\Delta \pi_{\alpha}(u(x,t)) \right]\pi_{\alpha}(u(x,t))\,\eta^2(x){\zeta(x,t)}\mu(x)\\
&\ge \frac12\sum_{x\in G}\,\rho(x)\,\partial_t \left[\pi_{\alpha}^2(u(x,t)){\zeta(x,t)}\right]\eta^2(x)\mu(x)\\
&\quad -\frac12\sum_{x\in G}\,\rho(x)\,\partial_t \zeta(x,t)\pi_{\alpha}^2(u(x,t))\,\eta^2(x)\mu(x)\\
&\quad -\frac1{4}\sum_{x,y\in G}\pi^2_{\alpha}(u(x,t))\left[1-\frac{\zeta(y,t)}{\zeta(x,t)}\right]^2 \eta^2(x){\zeta(x,t)}\,\omega(x,y) \\
&\quad - \sum_{x,y\in G}\pi_{\alpha}^2(u(x,t))[\eta(y)-\eta(x)]^2{\zeta(y,t)}\omega(x,y) \,.
\end{aligned}
\end{equation}

Moreover, observe that, due to the assumption $p\ge2$, $s\mapsto \pi_{\alpha}(s)$ defined in \eqref{eq43} is convex in $\R$ and of class $C^2$. Therefore, from \eqref{e5f} we get
\begin{equation}\label{eq417}
\Delta \pi_{\alpha}(u(x,t))\ge \pi_{\alpha}'(u(x,t))\Delta u(x,t)\quad \text{ for all }\; x\in G,\,\,t\in(0,T)\,.
\end{equation}
Thus, the left-hand-side of \eqref{eq416}, by \eqref{eq417} and \eqref{problema} satisfies the inequality
\begin{equation}\label{eq418}
\begin{aligned}
\sum_{x\in G}&\left[\rho(x)\partial_t \pi_{\alpha}(u(x,t))-\Delta \pi_{\alpha}(u(x,t)) \right]\pi_{\alpha}(u(x,t))\eta^2(x){\zeta(x,t)}\mu(x)\\
&\le \sum_{x\in G}\pi'_\alpha(u(x,t))\left[\rho(x) u_t(x,t)-\Delta u(x,t) \right]\pi_{\alpha}(u(x,t))\eta^2(x){\zeta(x,t)}\mu(x)= 0,
\end{aligned}
\end{equation}
for all $x\in G,\,t\in[0,T]$.
By combining together \eqref{eq416} and \eqref{eq418} we get
\begin{equation*}
\begin{aligned}
\frac{\partial}{\partial t}&\left(\sum_{x\in G}\,\rho(x)\,\pi_{\alpha}^2(u(x,t)){\zeta(x,t)}\eta^2(x)\mu(x)\right)\\
&\le \sum_{x\in G}\,\rho(x)\,\partial_t \zeta(x,t)\pi^2_{\alpha}(u(x,t))\,\eta^2(x)\mu(x)\\
&\quad +\frac1{2}\sum_{x,y\in G}\pi^2_{\alpha}(u(x,t))\left[1-\frac{\zeta(y,t)}{\zeta(x,t)}\right]^2 \eta^2(x){\zeta(x,t)}\,\omega(x,y) \\
&\quad + 2\sum_{x,y\in G}\pi_{\alpha}^2(u(x,t))\left[\eta(y)-\eta(x)\right]^2{\zeta(y,t)}\,\omega(x,y)\,.
\end{aligned}
\end{equation*}
We now integrate the latter with respect to the variable $t$ in the interval $[0,\tau]$, for any $\tau\in(0,T]$, thus we write
\begin{equation}\label{eq421}
\begin{aligned}
\sum_{x\in G}\,\rho(x)\,&\pi_{\alpha}^2(u(x,\tau)){\zeta(x,\tau)}\eta^2(x)\mu(x)-\sum_{x\in G}\,\rho(x)\alpha^{p/2}{\zeta(x,0)}\eta^2(x)\mu(x)\\
&\le\int_0^\tau \sum_{x\in G}\,\rho(x)\,\partial_t \zeta(x,t)\pi^2_{\alpha}(u(x,t))\,\eta^2(x)\mu(x)\,dt\\
&\quad +\frac1{2}\int_0^\tau\sum_{x,y\in G}\pi^2_{\alpha}(u(x,t))\left[1-\frac{\zeta(y,t)}{\zeta(x,t)}\right]^2 \eta^2(x){\zeta(x,t)}\,\omega(x,y)\,dt \\
&\quad + 2\int_0^\tau\sum_{x,y\in G}\pi_{\alpha}^2(u(x,t))\left[\eta(y)-\eta(x)\right]^2{\zeta(y,t)}\,\omega(x,y)\,dt\,,
\end{aligned}
\end{equation}
where we have used that $u(x,0)=0$ by \eqref{problema}. Now, since
\[\pi_\alpha(u)\to |u|^{\frac p2}, \quad \text{ as } \alpha\to 0^+,\]
letting $\alpha\to 0^+$ in \eqref{eq421}, we arrive to
\begin{equation*}\label{eq312}
\begin{aligned}
\sum_{x\in G}\,\rho(x)\,&|u(x,\tau)|^p\,{\zeta(x,\tau)}\eta^2(x)\mu(x)\\
&\le\int_0^\tau \sum_{x\in G}\,\left\{\rho(x)\,\partial_t \zeta(x,t)\mu(x)+\frac{1}2\zeta(x,t)\sum_{y\in G}\left[1-\frac{\zeta(y,t)}{\zeta(x,t)}\right]^2\omega(x,y)\right\}|u(x,t)|^p\,\eta^2(x)\,dt\\
&\quad + 2\int_0^\tau\sum_{x,y\in G}|u(x,t)|^p\left[\eta(y)-\eta(x)\right]^2{\zeta(y,t)}\,\omega(x,y)\,dt\,.
\end{aligned}
\end{equation*}
\end{proof}

\begin{proposition}\label{ineqv}
Let $(G, \mu, \omega)$ be a weighted graph such that \eqref{e7f} is satisfied. Let $u$ be a solution to problem \eqref{problema}. Let $p\ge1$ and $v\in\mathfrak F_T$ be a nonnegative function with finite support. Then, for any $\tau\in(0,T]$
\begin{equation}\label{eqlemma}
\int_0^\tau\sum_{x\in G}\left[\rho(x)\partial_t v(x,t)+\Delta v(x,t)\right]|u(x,t)|^p \mu(x)\,dt\ge\sum_{x\in G}\rho(x)|u(x,\tau)|^pv(x,\tau)\,\mu(x)
\end{equation}
\end{proposition}

\begin{proof}
For any $\alpha>0$, we define
\begin{equation}\label{G}
\pi_{\alpha}(u):=(u^2+\alpha)^{\frac p2}\,.
\end{equation}
By applying twice \eqref{e4f}, we obtain
\begin{equation}\label{eq33b}
\begin{aligned}
\sum_{x\in G}&\left[\rho(x)\,\partial_t \pi_\alpha(u(x,t))-\Delta\pi_\alpha(u(x,t))\right]v(x,t) \mu(x)\\
&=\sum_{x\in G}\left[-\rho(x)\,\partial_t v(x,t)-\Delta v(x,t)\right]\pi_\alpha(u(x,t)) \mu(x)\\&+\partial_t\left[\sum_{x\in G}\rho(x)\pi_\alpha(u(x,t))v(x,t)\,\mu(x)\right]\,.
\end{aligned}
\end{equation}
Since $p\geq 1$, the function $s\mapsto \pi_\alpha(s)$ in \eqref{G} is convex in $(0, +\infty)$ and of class $C^2$. Therefore, from \eqref{e5f} we get
\begin{equation}\label{e33f}
\Delta\pi_\alpha(u(x,t))\geq \pi'_\alpha(u(x,t))\Delta u(x,t)\quad \text{ for all }\; x\in G,\,t\in[0,T]\,.
\end{equation}
Moreover, due to \eqref{problema}
\begin{equation}\label{e33g}
\begin{aligned}
\rho(x)\partial_t \pi_\alpha(u(x,t))&=\pi'_\alpha(u(x,t))\,\rho(x) u_t(x,t)\\
&=\pi'_\alpha(u(x,t))\Delta u(x,t)\quad \text{for all}\,\, x\in G,\,t\in[0,T]\,.
\end{aligned}
\end{equation}
Therefore, by \eqref{eq33b}, \eqref{e33f} and \eqref{e33g} we obtain
\begin{equation*}
\sum_{x\in G}\left[\rho(x)\partial_t v(x,t)+\Delta v(x,t)\right]\pi_\alpha(u(x,t)) \mu(x)\ge\partial_t\left[\sum_{x\in G}\rho(x)\pi_\alpha(u(x,t))v(x,t)\,\mu(x)\right]
\end{equation*}
By integrate the latter inequality with respect to $t$ on $[0,\tau]$, for any $\tau\in (0,T]$, we get
\begin{equation}\label{e34f}
\begin{aligned}
\int_0^\tau\sum_{x\in G}&\left[\rho(x)\partial_t v(x,t)+\Delta v(x,t)\right]\pi_\alpha(u(x,t)) \mu(x)\,dt\\
&\ge\sum_{x\in G}\rho(x)\pi_\alpha(u(x,\tau))v(x,\tau)\,\mu(x)-\sum_{x\in G}\alpha^{p/2}\rho(x)v(x,0)\,\mu(x),
\end{aligned}
\end{equation}
where we have used that $u(x,0)=0$ by \eqref{problema}. Letting $\alpha\to 0^+$ in \eqref{e34f}, 
 $$
 \pi_\alpha(u(x,t))\to |u(x,t)|^p,
 $$
 hence we get the thesis.
 \end{proof}

\section{Some useful functions and their properties}\label{funzioni}\setcounter{equation}{0}

Let $x_0\in G$ be fixed. For any $R>0$ and $0<\delta<1/2$ we define the "cut-off" function
\begin{equation}\label{eq316}
\eta(x):=\min\left\{\frac{\left[R-s-d(x,x_0)\right]_+}{\delta R},\,\,1\right\}\, \quad \text{ for any }\; x\in G\,.
\end{equation}
The function $\eta$ will come into play in the proofs of all the theorems. We state in the following Lemma some of its properties. 
\begin{lemma}\label{lemma3}
The function $\eta$ defined in \eqref{eq316} satisfies
\begin{equation}\label{e55f}
|\dif \eta| \leq \frac 1{\delta R} d(x,y)\chi_{\{(1-\delta)R-2s\le d(x,x_0)\le R\}} \quad \text{ for any }\; x, y\in G,\,{ \omega(x,y)>0}\,.
\end{equation}
and
\begin{equation*}\label{eq317}
\sum_{y\in G} \left(\dif \eta\right)^2 \omega(x,y)\le \frac{1}{\delta^2 R^2}\mu(x)\chi_{\{(1-\delta)R-2s\le d(x,x_0)\le R\}} \quad \text{ for any }\; x\in G\,.
\end{equation*}
Furthermore,
\begin{equation}\label{e50f}
|\Delta \eta(x)|\leq \frac {C_0}{\delta R} \chi_{\{(1-\delta)R-2s\le d(x,x_0)\le R\}} \quad \text{ for any }\; x\in G\,.
\end{equation}
provided that {$d$} is also $1$-intrinsic with bound $C_0$ in the sense of Definition \ref{intrmet} with $q=1$.
\end{lemma}
The proof of Lemma \ref{lemma3} can be obtained with minor modifications to the proof of \cite[Lemma 2.3]{Huang}. 
\smallskip

\noindent The remainder of this section is devoted to introducing the test functions that we will use to prove our main results. These test functions are selected based on the class of uniqueness we are considering, which in turn depends on the assumptions made about $\rho$. Therefore, we need to consider separately the cases where the density satisfies either assumption \eqref{rhobound} or \eqref{rhobound3}.

\subsection{Density bounded away from zero}
 We define the functions
\begin{equation}\label{eq313}
\zeta(x,t)=e^{\xi(x,t)}\quad \text{ for all }\; x\in G,\,\, t\in[0,T]\,,
\end{equation}
and
\begin{equation}\label{xi}
\xi(x,t):=-\gamma t-\alpha\left[d(x,x_0)-\delta R\right]_+\quad \text{ for all }\; x\in G,\,\, t\in[0,T]\,,
\end{equation}
for some $\alpha,\gamma\ge0$, $R>0$, $0<\delta<\frac 12$. Regarding these two functions we show the following lemmas which will play a key role in the proof of Theorems \ref{teo2} and \ref{teo1} respectively.

\begin{lemma}\label{lemma41}
Let $p\geq 2$ and assume \eqref{rhobound}. Let $\xi$ be the function defined in \eqref{xi} with
$$
\alpha\ge0\quad\text{and}\quad\gamma\ge\frac{\alpha^2\,e^{2\alpha s}}{2\rho_0}\,,
$$
where $\rho_0$ has been defined in \eqref{rhobound}. Then
\begin{equation}\label{eq422}
\rho(x)\partial_t\xi(x,t)\,\mu(x)+\frac 12 \sum_{y\in G} \omega(x,y)\left[1-e^{\xi(y,t)-\xi(x,t)}\right]^2 \leq 0 \,\quad \text{ for all }\, x\in G,\,\, t\in[0,T]\,.
\end{equation}
\end{lemma}

\begin{lemma}\label{lemma2}
Let $p\geq 1$, $d$ be $1-$intrinsic as defined in Definition \ref{intrmet} and assume \eqref{rhobound}. Let $\zeta$ be the function defined in \eqref{eq313} with
$$
\alpha\ge0\quad\text{and}\quad \gamma\ge\frac{\alpha\,C_0\,e^{\alpha s}}{\rho_0}\,,
$$
where $\rho_0$ has been defined in \eqref{rhobound} and $C_0$ is the positive constant given in Definition \ref{intrmet}. Then
\begin{equation}\label{eq314}
\rho(x)\partial_t\zeta(x,t)+\Delta\zeta(x,t) \le 0 \quad \text{ for all }\, x\in G,\,\, t\in[0,T]\,.
\end{equation}
\end{lemma}

\begin{proof}[Proof of Lemma \ref{lemma41}]
It is direct to see that
\begin{equation}\label{eq423}
(e^a-1)^2\le e^{2|a|}a^2\quad \text{ for every }\; a\in \mathbb R\,.
\end{equation}
Moreover, by the triangle inequality for the distance $d$,
\[|\xi(x,t)-\xi(y,t)|\leq \alpha d(x,y) \quad \text{ for all } x,y\in G\,.\]
Therefore, by intrinsicity of $d$
\begin{equation}\label{eq424}
\begin{aligned}
\frac 12 \sum_{y\in G} \left[1-e^{\xi(y,t)-\xi(x,t)}\right]^2\omega(x,y)  &\le \frac 12 \sum_{y\in G}\left[\xi(y,t)-\xi(x,t)\right]^2e^{2|\xi(y,t)-\xi(x,t)|}\,\omega(x,y)\\
&\le \frac 12 \alpha^2\sum_{y\in G}d^2(x,y) e^{2\alpha d(x,y)}\omega(x,y)\\
&\le \frac {1}{2} \alpha^2e^{2s \alpha}\mu(x)\,,
\end{aligned}
\end{equation}
where we have used also the fact that
$
\omega(x,y)>0\iff x\sim y\,.
$
On the other hand, by the very definition of $\xi$ in \eqref{xi}, we can write
\begin{equation}\label{eq425}
\rho(x)\,\partial_t\xi(x,t)\,\mu(x)=-\gamma\rho(x)\mu(x)\,.
\end{equation}
By combining \eqref{eq424} and \eqref{eq425}, we get
$$
\begin{aligned}
\rho(x)\,\partial_t\xi(x,t)\,\mu(x)+\frac 12 \sum_{y\in G} \omega(x,y)\left[1-e^{\xi(y,t)-\xi(x,t)}\right]^2 &\leq -\gamma\rho(x)\mu(x) +\frac {1}{2} \alpha^2e^{2s \alpha}\mu(x)\\
&\le\left(-\gamma\rho_0 +\frac {1}{2} \alpha^2e^{2s \alpha}\right)\mu(x)\,,
\end{aligned}
$$
which ensure the thesis since $\gamma\ge\frac{\alpha^2\,e^{2\alpha s}}{2\rho_0}$.
\end{proof}

\begin{proof}[Proof of Lemma \ref{lemma2}]
In view of \eqref{eq313} and \eqref{rhobound}, we have
\begin{equation}\label{partialt}
\rho(x)\partial_t\zeta(x,t)=-\gamma\rho(x)\zeta(x,t)\le-\gamma\rho_0\zeta(x,t)\quad \text{ for any }\; x\in G,\,\, t\in [0,T]\,,
\end{equation}
and, due to \eqref{eq11}
\begin{equation}\label{eq420bis}
\begin{aligned}
\Delta \zeta(x,t)&=\frac1{\mu(x)}\sum_{y\in G}\omega(x,y)\left[\zeta(y,t)-\zeta(x,t) \right]\\
&\le\frac{1}{\mu(x)} \sum_{y\in G}\omega(x,y)|\nabla_{xy}\zeta| \quad \text{ for any }\; x\in G,\,\, t\in [0,T]\,.
\end{aligned}
\end{equation}
Let us now estimate 
$$
|\nabla_{xy}\zeta| =\left|e^{-\gamma t-\alpha [d(y, x_0)-\delta R]_+}-e^{-\gamma t-\alpha [d(x, x_0)-\delta R]_+}\right|
$$
distinguishing between the following four cases:
\begin{itemize}
\item[(a)] if $d(x,x_0)\le \delta R$ and $d(y,x_0)\le \delta R$ then $|\nabla_{xy}\zeta| =0$;
\item[(b)] if $d(x,x_0)\le \delta R$ and $d(y,x_0)> \delta R$ then, by \eqref{eq423}
$$
\begin{aligned}
|\nabla_{xy}\zeta| &=e^{-\gamma t}|e^{-\alpha [d(y, x_0)-\delta R]}-1|\\
&\le e^{-\gamma t}\alpha  [d(y, x_0)-\delta R]\,e^{\alpha [d(y, x_0)-\delta R]}\\
&\le e^{-\gamma t}\alpha [d(y, x_0)- d(x, x_0)]e^{\alpha [d(x, x_0)-\delta R]+\alpha [d(y, x_0)-d(x,x_0)]}\\
&\le \zeta(x,t)\alpha d(x, y)\,e^{\alpha d(x,y)}\,;
\end{aligned}
$$
\item[(c)] if $d(x,x_0)> \delta R$ and $d(y,x_0)\le \delta R$ then, by \eqref{eq423}
$$
\begin{aligned}
|\nabla_{xy}\zeta| &=e^{-\gamma t}|1-e^{-\alpha [d(x, x_0)-\delta R]}|\\
&\le e^{-\gamma t}\,e^{-\alpha [d(x, x_0)-\delta R]}|e^{\alpha [d(x, x_0)-\delta R]}-1|\\
&\le e^{-\gamma t-\alpha [d(x, x_0)-\delta R]_+} \alpha [ d(x, x_0)-\delta R]e^{\alpha [d(x, x_0)-\delta R]}\\
&\le \zeta(x,t) \alpha [d(x, x_0)- d(y, x_0)]e^{\alpha [d(x, x_0)-d(y,x_0)]}\\
&\le\zeta(x,t) \alpha d(x, y)\,e^{\alpha d(x,y)}\,;
\end{aligned}
$$
\item[(d)] if $d(x,x_0)> \delta R$ and $d(y,x_0)> \delta R$ then 
$$
\begin{aligned}
|\nabla_{xy}\zeta| &=e^{-\gamma t}|e^{-\alpha [d(y, x_0)-\delta R]}-e^{-\alpha [d(x, x_0)-\delta R]}|\\
&\le e^{-\gamma t -\alpha [d(x, x_0)-\delta R]}|e^{\alpha [d(x, x_0)-d(y,x_0)]}-1|\\
&\le \zeta(x,t) \alpha |d(x, x_0)- d(y, x_0)|\,e^{\alpha |d(x, x_0)-d(y,x_0)|}\\
&\le\zeta(x,t) \alpha d(x, y)\,e^{\alpha d(x,y)}\,.
\end{aligned}
$$
\end{itemize}
Implementing the latter cases into \eqref{eq420bis} and by \eqref{e14f}, \eqref{e1f} with $q=1$ we get
\begin{equation}\label{e43f}
\begin{aligned}
\Delta \zeta(x,t)&=\frac1{\mu(x)}\sum_{y\in G}\omega(x,y)\left[\zeta(y,t)-\zeta(x,t) \right]\\
&\le\alpha\frac{\zeta(x,t)}{\mu(x)} \sum_{y\in G}\omega(x,y) d(x, y)\,e^{\alpha d(x,y)}\\
&\le \alpha\,C_0\,e^{\alpha s} \zeta(x,t) 
\end{aligned}
\end{equation}
From \eqref{partialt} and \eqref{e43f}, we get for any $x\in G$, $t\in[0,T]$
\begin{equation*}
\rho(x)\partial_t\zeta(x,t)+\Delta\zeta(x,t) 
\le\left(-\gamma\rho_0+ \alpha\,C_0\,e^{\alpha s} \right)\zeta(x,t)\,.
\end{equation*}
The thesis follows because $\gamma\ge\frac{\alpha\,C_0\,e^{\alpha s}}{\rho_0}$ and $\zeta(x,t)>0$ for any $x\in G$ and $t\ge0$.
\end{proof}

\subsection{Vanishing density}
In this second part of our investigation we deal with the following test function
\begin{equation}\label{eq51}
\theta(x,t)=e^{-\gamma t}[d(x,x_0)+k]^{-\alpha}\quad \text{ for all }\; x\in G,\,\, t\in[0,T]\,,
\end{equation}
where $\alpha,\gamma\ge0$ and $k>\max\{1,s\}$ is the same of assumption \eqref{rhobound3}. Similarly to the previous subsection, also in this case it is useful to show that $\theta$ satisfies some inequalities. For this reason we state the following two Lemmas which will be implemented in the proofs of Theorem \ref{teo3} and \ref{teo4}, respectively.
\begin{lemma}\label{lemma62}
Let $(G, \mu, \omega)$ be a weighted graph as defined in \eqref{e7f}, let $p \geq 2$ and assume \eqref{rhobound3} with $0<\sigma\le2$. Let $\theta$ be the function defined in \eqref{eq51} with
$$
\alpha\ge0\quad\text{and}\quad \gamma\ge\frac{\alpha^2\,C_1^2}{2\rho_0}\,.
$$
where $\rho_0$ has been defined in \eqref{rhobound3} and 
 \begin{equation}\label{eq64}
C_1 := \frac{k^{\alpha + 1}}{(k-s)^{\alpha + 1}}\,.
\end{equation}Then
\begin{equation}\label{eq62}
\begin{aligned}
 \rho(x)&\partial_t\theta(x,t)\mu(x)+\frac{\theta(x,t)}{2}\sum_{y \in G}\omega(x,y) \left[1- \frac{\theta(y,t)}{\theta(x,t)}\right]^2 \le 0\quad \text{for all}\,\,x\in G, \,\,t\in[0,T].
 \end{aligned}
\end{equation}
\end{lemma}

\begin{lemma}\label{lemma51}
Let $(G, \mu, \omega)$ be a weighted graph as defined in \eqref{e7f}, let $p\geq 1$, $d$ be $1-$intrinsic as defined in Definition \ref{intrmet} and assume \eqref{rhobound3} with $0<\sigma\le1$. Let $\theta$ be the function defined in \eqref{eq51} with
$$
\alpha\ge0\quad\text{and}\quad \gamma\ge\frac{\alpha\,C_0\,C_1}{\rho_0}\,.
$$
where $\rho_0$ has been defined in \eqref{rhobound3}, $C_0$ in Definition \ref{intrmet} and $C_1$ in \eqref{eq64}. Then
\begin{equation}\label{eq52}
\rho(x)\partial_t\theta(x,t)+\Delta\theta(x,t) \le 0 \quad \text{for all}\,\,\, x\in G,\,\, t\in[0,T].
\end{equation}
\end{lemma}

Before proving the two lemmas, we state the following result which has been showed in \cite[Lemma 4.1]{Svagna}.
\begin{lemma}\label{lemma61}
Let $(G,\omega,\mu)$ be as in \eqref{e7f}, let $x,y \in G$ with $\omega(x,y) > 0$ be given and $\alpha>0$. Then
\begin{equation}\label{eq61}
\big | [d(y,x_0) + k]^{-\alpha} - [d(x,x_0) + k]^{-\alpha} \big | \leq \alpha d(x,y) [d(x,x_0) + k - s]^{-\alpha -1},
\end{equation}
for some $k>\max\{1,s\}$.
\end{lemma}

\begin{proof}[Proof of Lemma \ref{lemma62}]
In view of the very definition of $\theta$ given in \eqref{eq51} and due to Lemma \ref{lemma61} we deduce that
\begin{equation}\label{eq67}
\begin{aligned}
    \sum_{y \in G}\omega(x,y) \left[ 1- \frac{\theta(y,t)}{\theta(x,t)}\right]^2& =\sum_{y \in G}\omega(x,y) \left(\frac{\left[(d(x,x_0)+k)^{-\alpha} - (d(y,x_0)+k)^{-\alpha} \right]}{(d(x,x_0)+k)^{-\alpha}}  \right)^2\\
    & \leq \sum_{y \in G}\omega(x,y) \frac{\left[\alpha d(x,y)(d(x,x_0) + k - s)^{-\alpha -1}\right]^2}{(d(x,x_0)+k)^{-2 \alpha}} \\
    &=\frac{\alpha^2(d(x,x_0) + k - s)^{-2\alpha -2}}{(d(x,x_0)+k)^{-2 \alpha}} \sum_{y \in G}\omega(x,y) d^2(x,y).
\end{aligned}
\end{equation}
Now, we want to find a constant $C_1$ such that
\begin{equation}\label{eq68}
    [d(x,x_0) + k-s]^{-\alpha-1} \leq C_1 [d(x,x_0) + k]^{-\alpha-1} \quad\text{for any}\,\,\, x  \in G.
\end{equation}
We therefore define the function $g :[0,+\infty) \longrightarrow \mathbb{R}$,
\begin{equation*}
    g(z):= \frac{[z+k-s]^{-\alpha-1}}{[z+k]^{-\alpha-1}}.
\end{equation*}
Then
$$
\begin{aligned}
   g'(z) &= \frac{- (\alpha +1)[z + k-s]^{-\alpha-2}[z+ k]^{-\alpha-1}+(\alpha +1) [z + k]^{-\alpha-2} [z + k-s]^{-\alpha-1} }{[z + k]^{-2(\alpha +1)}}\\
   & =\frac{-s(\alpha +1)[z + k-s]^{-\alpha-2} }{[z + k]^{-\alpha}}\le 0 \quad\text{for all}\,\,\,z\in[0,+\infty)\,.
\end{aligned}
$$
Therefore the function $g(z)$ is strictly decreasing and attains its maximum in $z= 0$. Consequently we may say that
\begin{equation}\label{eq66}
\frac{[d(x,x_0) + k-s]^{-\alpha-1}}{[d(x,x_0) + k]^{-\alpha-1}}=g(d(x,x_0))\le g(0)= \frac{(k-s)^{-\alpha-1}}{k^{-\alpha-1}}
\end{equation}
By virtue of \eqref{eq66}, \eqref{eq68} and due to the intrinsicity property of the metric $d$, inequality \eqref{eq67}, by also using \eqref{eq64} reduces to
\begin{equation}\label{eq65}
\begin{aligned}
    \sum_{y \in G}\omega(x,y) \left[ 1- \frac{\theta(y,t)}{\theta(x,t)}\right]^2 
    & \leq\frac{\alpha^2(d(x,x_0) + k )^{-2\alpha -2}}{(d(x,x_0)+k)^{-2 \alpha}} \sum_{y \in G}\omega(x,y) d^2(x,y).\\
& \leq C_1^2 \alpha ^2 \big[(d(x,x_0) + k)]^{-2} \mu(x).
\end{aligned}
\end{equation}
We finally estimate the sum
$$
 \rho(x)\partial_t\theta(x,t)+\frac{\theta(x,t)}{2}\sum_{y \in G}\omega(x,y) \left[ 1- \frac{\theta(y,t)}{\theta(x,t)}\right]^2\,.
$$
Due to the definition of $\theta$ in \eqref{eq51}, \eqref{eq65} and \eqref{rhobound3}, the latter reduces to
$$
\begin{aligned}
\rho(x)\partial_t\theta(x,t)\mu(x)&+\frac{\theta(x,t)}{2}\sum_{y \in G}\omega(x,y) \left[ 1- \frac{\theta(y,t)}{\theta(x,t)}\right]^2\\
&\le -\gamma\rho(x)e^{-\gamma t}[d(x,x_0)+k]^{-\alpha}\mu(x)+\frac12C_1^2 \alpha ^2 e^{-\gamma t}[d(x,x_0) + k]^{-2-\alpha} \mu(x)\\
&\le \theta(x,t)\mu(x)\left\{-\rho_0\gamma[d(x,x_0)+k]^{-\sigma}+\frac12C_1^2 \alpha ^2 \big[d(x,x_0) + k]^{-2}\right\}\,.
\end{aligned}
$$
In view of the assumptions on $\sigma$ and $\gamma$ we get the thesis.
\end{proof}

\begin{proof}[Proof of Lemma \ref{lemma51}]
In view of \eqref{eq51} and \eqref{rhobound3}, we have
\begin{equation}\label{eq53}
\begin{aligned}
\rho(x)\partial_t\theta(x,t)&=-\gamma\rho(x)\theta(x,t)\\
&\le-\gamma\rho_0[d(x,x_0)+k]^{-\sigma}\theta(x,t)\quad \text{for any}\,\,x\in G,\,\, t\in [0,T],\end{aligned}
\end{equation}
and
\begin{equation*}
\begin{aligned}
\Delta \theta(x,t)&=\frac1{\mu(x)}\sum_{y\in G}\omega(x,y)\left[\theta(y,t)-\theta(x,t) \right]\\
&\le\frac{1}{\mu(x)} \sum_{y\in G}\omega(x,y)|\nabla_{xy}\theta| \quad \text{for any}\,\, x\in G,\,\, t\in [0,T]\,.
\end{aligned}
\end{equation*}
We now estimate 
$$
|\nabla_{xy}\theta| = e^{-\gamma t}\left|[d(y, x_0)+k]^{-\alpha}-[d(x, x_0)+k]^{-\alpha}\right|\,.
$$
By means of Lemma \ref{lemma61}, we get
\begin{equation}\label{eq55}
\begin{aligned}
     |\nabla_{xy}\theta| &\le e^{-\gamma t} \alpha [d(x,x_0) + k-s]^{-\alpha-1} d(x,y) \\
\end{aligned}
\end{equation}
By combining \eqref{eq53} and \eqref{eq55}, together with the is $1-$intrinsic assumption on the pseudo metric $d$, we get
\begin{equation}\label{eq56}
\begin{aligned}
 &\rho(x)\partial_t\theta(x,t)+\Delta\theta(x,t) \\
 &\leq -\gamma\rho_0[d(x,x_0)+k]^{-\sigma}\theta(x,t) + e^{-\gamma t} \alpha   [d(x,x_0) + k-s]^{-\alpha-1}\frac{1}{\mu(x)}\sum_{y\in G}\omega(x,y) d(x,y)\\
 &\le -\gamma\rho_0[d(x,x_0)+k]^{-\sigma}\theta(x,t) + e^{-\gamma t} \alpha C_0  [d(x,x_0) + k-s]^{-\alpha-1}
\end{aligned}
\end{equation}
In view of \eqref{eq68} with $C_1$ as in \eqref{eq64}, we get
\begin{equation}
\begin{aligned}
 \rho(x)\partial_t\theta(x,t)&+\Delta\theta(x,t) \\
 & \leq -\gamma\rho_0[d(x,x_0)+k]^{-\sigma}\theta(x,t) + e^{-\gamma t} \alpha C_0C_1  [d(x,x_0) + k]^{-\alpha-1}\\
 &= \theta(x,t)\left[-\gamma\rho_0[d(x,x_0)+k]^{-\sigma} +  \alpha C_0C_1  [d(x,x_0) + k]^{-1}\right]\\
\end{aligned}
\end{equation}
The thesis follows by assumptions on $\sigma$ and $\gamma$.
\end{proof}

\section{Proof of Theorems \ref{teo2} and \ref{teo1}}\label{sec5}\setcounter{equation}{0}

We are ready to prove the results of Subsection \ref{bounded}. We recall the reader that here $\rho$ satisfies assumption \eqref{rhobound}, i.e. it is a positive function bounded from below by the positive constant $\rho_0$.
 
\begin{proof}[Proof of Theorem \ref{teo2}]
Let $\eta$ and $\xi$ be defined by  \eqref{eq316} and \eqref{xi} respectively, 
with
\begin{equation}\label{coeffassumption2}
\alpha=\frac{K}{2s}\log(R_n),\quad \text{and}\quad \gamma:=\frac{\alpha^2 e^{2\alpha s}}{2\rho_0},
\end{equation}
where $K\in\left(\frac{2\beta s}{1-2\delta},1\right)$, $\delta\in\left(0,\frac12-\beta s\right)$ and $\rho_0$ has been defined in \eqref{rhobound}. Observe that $\left(0,\frac12-\beta s\right)$ is non empty due to the assumption on $\beta$ in Theorem \ref{teo2}. Moreover, chose $R>\frac{2s}{1-2\delta}$ in \eqref{eq316} and \eqref{xi}.
From Proposition \ref{prop41} with $\eta$ as above and $\zeta(x,t)=e^{\xi(x,t)}$ for any $(x,t)\in S_T$ we get
$$
\begin{aligned}
\sum_{x\in G}&|u(x,\tau)|^p\rho(x)\eta^2(x)e^{\xi(x,\tau)}\mu(x)\le2\int_0^\tau \sum_{x,y\in G}|u(x,t)|^pe^{\xi(y,t)}\left[\eta(y)-\eta(x)\right]^2\omega(x,y)\,dt\\
&+\int_0^\tau \sum_{x\in G}|u(x,t)|^p\eta^2(x)e^{\xi(x,t)}\left\{\rho(x)\partial_t\xi(x,t)\mu(x)+\frac 12\sum_{y\in G}\omega(x,y)\left[1-e^{\xi(y,t)-\xi(x,t)}\right]^2 \right\}\,dt\,.
\end{aligned}
$$
By applying Lemma \ref{lemma41}, we get rid of the second term in the right-hand-side, thus we write
\begin{equation}\label{eq426}
\sum_{x\in G}|u(x,\tau)|^p\rho(x)\eta^2(x)e^{\xi(x,\tau)}\mu(x)\le2\int_0^\tau \sum_{x,y\in G}|u(x,t)|^pe^{\xi(y,t)}\left[\eta(y)-\eta(x)\right]^2\omega(x,y)\,dt\,.
\end{equation}
Since $R>\frac{2s}{1-2\delta}$, we have $\eta(x)=1$ when $d(x,x_0)\le\delta R$. Therefore, for the left hand side we write
\begin{equation}\label{eq427}
e^{-\gamma \tau} \rho_0\sum_{x\in B_{\delta R}(x_0)}|u(x,\tau)|^p\,\mu(x)\le  \sum_{x\in G}\rho(x)|u(x,\tau)|^p\eta(x)e^{\xi(x,\tau)}\,\mu(x)\,.
\end{equation} 
For the right-hand-side of \eqref{eq426} we first observe that
$$
\xi(y,t)\le\alpha\,s\,+\xi(x,t)\,
$$
with $s$ being the jump-size defined in \eqref{e14f}. Therefore, by Lemma \ref{lemma2} and since $d$ is intrinsic we get
\begin{equation}\label{eq428}
\begin{aligned}
\sum_{x,y\in G}&|u(x,t)|^pe^{\xi(y,t)}\left[\eta(y)-\eta(x)\right]^2\omega(x,y)\le e^{\alpha s}\sum_{x,y\in G}|u(x,t)|^pe^{\xi(x,t)}|\nabla_{xy}\eta|^2\omega(x,y)\\
&\le \frac{e^{\alpha s}}{(\delta R)^2}\sum_{x\in G}|u(x,t)|^pe^{-\alpha[d(x,x_0)-\delta R]}\chi_{\{(1-\delta)R-2s\le d(x,x_0)\le R\}}\sum_{y\in G}d^2(x,y)\omega(x,y)\\
&\le \frac{e^{\alpha s}}{(\delta R)^2}e^{\alpha\delta R-\alpha[(1-\delta)R-2s]}\sum_{x\in B_R}|u(x,t)|^p\mu(x)\,.
\end{aligned}
\end{equation}
We substitute \eqref{eq427} and \eqref{eq428} into \eqref{eq426}, thus we have
\begin{equation}\label{eq426b}
 \sum_{x\in B_{\delta R}(x_0)}|u(x,\tau)|^p\,\mu(x)\le \frac{2}{\rho_0(\delta R)^2}e^{\gamma \tau+3\alpha s-(1-2\delta) \alpha R}\int_0^\tau\sum_{x\in B_R}|u(x,t)|^p\mu(x)\,dt.
\end{equation}
Let $R_n$ be as in assumption (iii) of Theorem \ref{teo2}. For $n$ big enough such that $R_n>\frac{2s}{1-2\delta}$, let $R=R_n$ in the latter inequality. Then, by using \eqref{maincondition2} we get
\begin{equation}\label{eq429}
 \sum_{x\in B_{\delta R_n}(x_0)}|u(x,\tau)|^p\,\mu(x)\le \frac{2}{\rho_0(\delta R_n)^2}\,\operatorname{exp}\{\gamma \tau+3\alpha s-(1-2\delta) \alpha R_n+\beta R_n\log R_n\}\,.
\end{equation}
Due to \eqref{coeffassumption2}, the exponent of the exponential reads
$$
\frac{1}{\rho_0}\frac{K^2}{8s^2}\,\tau 
R_n^K\log^2 R_n+\frac 32 K\,\log R_n-\left[(1-2\delta)\frac K{2s}-\beta\right]R_n\log R_n=:f(R_n).
$$
We now take the limit as $n\to+\infty$ in \eqref{eq429}. Since $K\in\left(\frac{2\beta s}{1-2\delta},1\right)$, $\delta\in\left(0,\frac12-\beta s\right)$, $\beta\in\left(0,\frac1{2s}\right)$ and $R_n\to+\infty$ as $n\to+\infty$, the right-hand-side of \eqref{eq429} gives
$$
\operatorname{exp}\{f(R_n)\}\to 0\quad \text{as}\,\,\,n\to+\infty\,.
$$
Therefore we get
$$
 \sum_{x\in G}|u(x,\tau)|^p\,\mu(x)\le0\,,
$$
which ensures that $u(x,\tau)=0$ for all $x\in G$. Note that $\tau$ is arbitrarily chosen in $(0,T]$ hence the theorem follows.

\end{proof}

\begin{proof}[Proof of Theorem \ref{teo1}]
Let $\eta$ and $\zeta$ be defined by  \eqref{eq316} and \eqref{eq313}, respectively with
\begin{equation}\label{coeffassumption}
\alpha=\frac{K}{s}\log R_n\quad\text{and}\quad \gamma:=\alpha e^{\alpha s}\frac{C_0}{\rho_0},
\end{equation}
where $\rho_0$ has been defined in \eqref{rhobound}.
Moreover, let
\begin{equation}\label{constTeo32}
K\in\left(\frac{\beta s}{1-2\delta},1\right),\quad 0<\delta<\frac{1-\beta s}{2},
\end{equation}
and choose $R>\frac{2s}{1-2\delta}$. Observe that $\left(0,\frac{1-\beta s}{2}\right)$ is non empty due to the assumption on $\beta$ made in Theorem \ref{teo1}. From Proposition \ref{ineqv} with $$v(x,t):=\eta(x)\zeta(x,t)\quad\text{for all}\,\,x\in G,\,\,t\in[0,T],$$
 we can write
 $$
\int_0^\tau\sum_{x\in G}\left[\rho(x)\eta(x)\partial_t \zeta(x,t)+\Delta (\eta(x)\zeta(x,t))\right]|u(x,t)|^p \mu(x)\,dt\ge\sum_{x\in G}\rho(x)|u(x,\tau)|^p\eta(x)\zeta(x,\tau)\,\mu(x).
 $$
Observe that our choice of $v$ has finite support in $G$ due to the presence of $\eta$. By using Lemma \ref{lemma5} we compute
 $$
 \Delta (\eta(x)\zeta(x,t))=\eta(x)\Delta\zeta(x,t)+\zeta(x,t)\Delta\eta(x)+\frac1{\mu(x)}\sum_{y\in G}\nabla_{x,y}\zeta\nabla_{xy}\eta\omega(x,y)\,.
 $$
 Therefore we get, by applying also Lemma \ref{lemma2}
 \begin{equation}\label{eq71}
 \begin{aligned}
 \sum_{x\in G}\rho(x)&|u(x,\tau)|^p\eta(x)\zeta(x,\tau)\,\mu(x)\\
 &\le \int_0^\tau\sum_{x\in G}\left[\rho(x)\partial_t \zeta(x,t)+\Delta \zeta(x,t)\right]|u(x,t)|^p \eta(x)\mu(x)\,dt+\int_0^\tau (I_1+I_2)\,dt\\
 &\le \int_0^\tau (I_1+I_2)\,dt
 \end{aligned}
 \end{equation}
where
\begin{align}
&I_1:=\sum_{x\in G}|u(x,t)|^p\zeta(x,t)\,\Delta\eta(x)\,\mu(x)\,,\label{eq34a}\\
&I_2:=\sum_{x,y\in G}|u(x,t)|^p\nabla_{xy}\zeta\,\nabla_{xy}\eta\,\omega(x,y)\,.\label{eq34b}
\end{align}
Note that since $R>\frac{2s}{1-2\delta}$, we have $\eta(x)=1$ when $d(x,x_0)\le\delta R$. Therefore, for the left hand side of \eqref{eq71}, we write
\begin{equation}\label{lhs71}
e^{-\gamma \tau} \rho_0\sum_{x\in B_{\delta R}(x_0)}|u(x,\tau)|^p\,\mu(x)\le  \sum_{x\in G}\rho(x)|u(x,\tau)|^p\eta(x)\zeta(x,\tau)\,\mu(x)\,.
\end{equation} 
We now consider the right hand side of \eqref{eq71}. By exploiting the definition of $\eta$ given in \eqref{eq316}, by \eqref{e50f} and assumption \eqref{e1f}, we get
\begin{equation}\label{eq56g}
\begin{aligned}
\left| \int_0^\tau I_{1}\,dt \right| &\le\int_0^\tau\sum_{x\in G}|u(x,t)|^p \zeta(x,t)|\Delta \eta(x)|\mu(x)\,dt\\
&\le \frac{C_0}{\delta R}\int_0^\tau\sum_{x\in G}|u(x,t)|^p \zeta(x,t)\mu(x)\chi_{\{(1-\delta)R-2s\le d(x,x_0)\le R\}}\,dt\\
&\le \frac{C_0}{\delta R}\int_0^\tau\sum_{x\in G}|u(x,t)|^p e^{-\alpha[d(x,x_0)-\delta R]_+}\mu(x)\chi_{\{(1-\delta)R-2s\le d(x,x_0)\le R\}}\,dt\\
&\le \frac{C_0}{\delta R}e^{2\alpha s-(1-2\delta)\alpha R}\int_0^\tau\sum_{x\in B_R(x_0)}|u(x,t)|^p \mu(x)\,dt\,.
\end{aligned}
\end{equation}
Let us now consider $I_2$. Arguing as in the proof of Lemma \ref{lemma2}, we obtain that
$$
|\nabla_{xy}\zeta|\le \zeta(x,t)\alpha\,d(x,y)\,e^{\alpha\,d(x,y)}\quad\text{for any}\,\,x,y\in G,\,\,t\in[0,T].
$$
Therefore, combining the latter with \eqref{e55f}, \eqref{e14f} and \eqref{e1f} we obtain
\begin{equation}\label{e47f}
\begin{aligned}
\left|\int_0^\tau I_{2}\,dt\right| &\le \int_0^\tau \sum_{x, y\in G} |u(x,t)|^p |\nabla_{xy}\zeta| |\nabla_{xy}\eta| \omega(x,y)\\
&\leq \frac{\alpha}{\delta R} \int_0^\tau \sum_{x,y\in G} |u(x,t)|^p \zeta(x,t) d^2(x,y)\,e^{\alpha\,d(x,y)}\chi_{\{(1-\delta)R-2s\le d(x,x_0)\le R\}}\omega(x,y) \\
&\leq \frac{ \alpha e^{\alpha s}}{\delta R}  \int_0^\tau  \sum_{x\in G} \zeta(x,t)|u(x,t)|^p \chi_{\{(1-\delta)R-2s\le d(x,x_0)\le R\}} \sum_{y\in G} d^2(x,y)\omega(x,y)\\
&\leq \frac{  C_0\alpha e^{\alpha s}}{\delta R}  \int_0^\tau  \sum_{x\in G} e^{-\alpha[d(x,x_0)-\delta R]_+}|u(x,t)|^p \chi_{\{(1-\delta)R-2s\le d(x,x_0)\le R\}} \mu(x)\\
&\leq \frac{ C_0 \alpha e^{3\alpha s-(1-2\delta)\alpha R}}{\delta R}  \int_0^\tau  \sum_{x\in B_R(x_0)} |u(x,t)|^p  \mu(x).
\end{aligned}
\end{equation}
Let $R_n$ be as in Theorem \ref{teo1}. We chose $R=R_n$; by condition \eqref{maincondition2},
$$
\int_0^\tau\sum_{x\in B_{R_n}(x_0)}|u(x,t)|^p \mu(x)\,dt\le C\,e^{\beta R_n\log R_n}.
$$
Hence, \eqref{eq56g} and \eqref{e47f} yield
\begin{equation}\label{rhs71}
\begin{aligned}
&\left| \int_0^\tau I_{1}\,dt \right|\le \frac{C_0C}{\delta R_n}e^{2\alpha s-(1-2\delta)\alpha R_n+\beta R_n\log R_n}\\
&\left| \int_0^\tau I_{2}\,dt \right|\le \frac{C_0C\alpha }{\delta R_n}e^{3\alpha s-(1-2\delta)\alpha R_n+\beta R_n\log R_n}
\end{aligned}
\end{equation}
Combining together \eqref{eq71}, \eqref{lhs71} and \eqref{rhs71} we get
\begin{equation}\label{last}
\rho_0 \sum_{x\in B_{\delta R_n}(x_0)}|u(x,\tau)|^p\,\mu(x)\le \frac{C_0C}{\delta R_n}\left[1+\alpha e^{\alpha s}\right]e^{\gamma \tau+2\alpha s-(1-2\delta)\alpha R_n+\beta R_n\log R_n}
\end{equation}
Due to \eqref{coeffassumption}, the exponent of the exponential reads
$$
\frac{C_0}{\rho_0}\frac{K}{s}\,\tau 
R_n^K\log R_n+2 K\,\log R_n-\left[(1-2\delta)\frac K{s}-\beta\right]R_n\log R_n=:g(R_n).
$$
We now take the limit as $n\to+\infty$ in \eqref{last}. Due to the choices of $K$, $\delta$ in \eqref{constTeo32}, 
since $\beta\in\left(0,\frac1{s}\right)$ and $R_n\to+\infty$ as $n\to+\infty$, we get
$$
\operatorname{exp}\{g(R_n)\}\to 0\quad \text{as}\,\,\,n\to+\infty\,.
$$
Therefore we get
$$
 \sum_{x\in G}|u(x,\tau)|^p\,\mu(x)\le0\,,
$$
which ensures that $u(x,\tau)=0$ for all $x\in G$. Note that $\tau$ is arbitrarily chosen in $(0,T]$ hence the theorem follows.
\end{proof}

\section{Proof of Theorems \ref{teo3} and \ref{teo4}}\label{sec6}\setcounter{equation}{0}

We are left to show the validity of the statements in Subsection \ref{vanishing}. We recall the reader that here $\rho$ satisfies assumption \eqref{rhobound3}, i.e. it is a positive function which may vanish at infinity according to the assumptions of each theorem.

\begin{proof}[Proof of Theorem \ref{teo3}]
Let $\eta$ and $\theta$ be defined by  \eqref{eq316} and \eqref{eq51} respectively, with
\begin{equation}\label{coeffassumption3}
\alpha >\max\{0,\beta-2\},\quad \text{and}\quad \gamma\ge\frac{(C_1 \alpha)^2}{2\rho_0},
\end{equation}
where $\rho_0$ has been defined in \eqref{rhobound3} and $C_1$ in \eqref{eq64}. 
From Proposition \ref{prop41} with $\eta$ as above and $\zeta(x,t)=\theta(x,t)$ for any $(x,t)\in S_T$ we get
$$
\begin{aligned}
\sum_{x\in G}&|u(x,\tau)|^p\rho(x)\eta^2(x){\theta(x,\tau)}\mu(x)\le2\int_0^\tau \sum_{x,y\in G}|u(x,t)|^p{\theta(y,t)}\left[\eta(y)-\eta(x)\right]^2\omega(x,y)\,dt\\
&+\int_0^\tau \sum_{x\in G}|u(x,t)|^p\eta^2(x)\left\{\rho(x)\partial_t\theta(x,t)\mu(x)+\frac {{\theta(x,t)}}2\sum_{y\in G}\omega(x,y)\left[1-\frac{\theta(y,t)}{\theta(x,t)}\right]^2 \right\}\,dt\,.
\end{aligned}
$$
By applying Lemma \ref{lemma62}, we get
\begin{equation}\label{eq69}
\sum_{x\in G}|u(x,\tau)|^p\rho(x)\eta^2(x)\theta(x,\tau)\mu(x)\le2\int_0^\tau \sum_{x,y\in G}|u(x,t)|^p{\theta(y,t)}\left[\eta(y)-\eta(x)\right]^2\omega(x,y)dt.
\end{equation}
Consider the right-hand-side of \eqref{eq69}, we observe that, by triangular inequality, $$d(y,x_0)\ge d(x,x_0)-s\quad\text{for any}\,\,\, x\sim y,$$ hence we can write
$$[d(y,x_0) + k]^{-\alpha} \leq [d(x,x_0) + k-s]^{-\alpha} \quad\text{for any}\,\,\,x \sim y.$$
The latter, together with the definition of $\theta$ in \eqref{eq51} yield,
\begin{equation}\label{eq610}
\begin{aligned}
 \sum_{x,y \in G} |u(x,t)|^p &\theta(y,t) [\eta(y)-\eta(x)]^2\omega(x,y) \\
 &= e^{-\gamma t}\sum_{x,y \in G} |u(x,t)|^p [d(y,x_0) + k]^{-\alpha}  [\eta(y)-\eta(x)]^2\omega(x,y)\\
 &\leq e^{-\gamma t}\sum_{x,y \in G} |u(x,t)|^p [d(x,x_0) + k-s]^{-\alpha}  [\eta(y)-\eta(x)]^2\omega(x,y).
\end{aligned}
\end{equation}
Moreover, we observe that arguing as in the proof of Lemma \ref{lemma62}, we can claim that there exists $C_2>0$ independently of $R$ such that 
\begin{equation}\label{eq69bis}
[d(x,x_0) + k-s]^{-\alpha} \leq C_2 \,[d(x,x_0) + k]^{-\alpha} \quad \text{ for any } x\in G.
\end{equation}
If we now apply Lemma \ref{lemma3} and the intrinsicity of $d$ to inequality \eqref{eq610}, for $C_2>0$ as in \eqref{eq69bis}, we have
\begin{equation}\label{eq611}
\begin{aligned}
& \sum_{x,y \in G} |u(x,t)|^p \theta(y,t) [\eta(y)-\eta(x)]^2\omega(x,y) \\
 &\le \frac{C_2}{\delta^2 R^2} e^{-\gamma t}\sum_{x\in G} |u(x,t)|^p [d(x,x_0) + k]^{-\alpha} \chi_{\{(1-\delta)R-2s \leq d(x,x_0)\leq R \} }\sum_{y\in G}d^2(x,y)\omega(x,y)\\
 &\le\frac{C_2}{\delta^2 R^2} [(1-\delta)R-2s+ k]^{-\alpha}e^{-\gamma t}\sum_{x\in B_R} |u(x,t)|^p  \mu(x)\,.
 \end{aligned}
\end{equation}
We now substitute \eqref{eq611} into \eqref{eq69}, thus we have
$$
\begin{aligned}
\sum_{x\in G}|u(x,\tau)|^p\rho(x)\eta^2(x)\theta(x,\tau)\mu(x)
&\le 
\frac{C_2}{\delta^2 R^2}[(1-\delta)R-2s+ k]^{-\alpha}\int_0^\tau e^{-\gamma t}\sum_{x\in B_R} |u(x,t)|^p \mu(x)\,dt\\
&\le 
\frac{C_2}{\delta^2 R^2}[(1-\delta)R-2s+ k]^{-\alpha}\int_0^\tau\sum_{x\in B_R} |u(x,t)|^p \mu(x)\,dt\\
\end{aligned}
$$
Let $R_n$ be as in the assumptions of Theorem \ref{teo3}. 
Then, by using \eqref{maincondition3} we get
\begin{equation}\label{eq612}
\begin{aligned}
\sum_{x\in G}|u(x,\tau)|^p\rho(x)\eta^2(x)\theta(x,\tau)\mu(x)
&\le\frac{C_2}{\delta^2 R_n^2}[(1-\delta)R_n-2s+ k]^{-\alpha}[R_n+k]^\beta\,.\end{aligned}
\end{equation}
We now take the limit as $n\to+\infty$ in \eqref{eq612}. Since $\alpha>\max\{0,\beta-2\}$, $R_n\to+\infty$ as $n\to+\infty$ and $\eta\to1$ as $R_n\to\infty$, the latter gives
$$
 \sum_{x\in G}|u(x,\tau)|^p\rho(x)\theta(x,\tau)\,\mu(x)\le0\,.
$$
Due to \eqref{rhobound3}, since $\theta>0$ in $S_T$ and $|u|^p\ge0$, the latter yield that $u(x,\tau)=0$ for all $x\in G$. Note that $\tau$ is arbitrarily chosen in $(0,T]$ hence the theorem follows.

\end{proof}


 \begin{proof}[Proof of Theorem \ref{teo4}]
Let $\eta$ and $\theta$ be defined by  \eqref{eq316} and \eqref{eq51}, respectively with
\begin{equation}\label{coeffassumption4}
\alpha>\max\{0,\beta-1\} \quad\text{and}\quad \gamma\ge\frac{\alpha\,C_0\,C_1}{\rho_0},
\end{equation}
where $\rho_0$ has been defined in \eqref{rhobound3}.
From Proposition \ref{ineqv} with $$v(x,t):=\eta(x)\theta(x,t)\quad\text{for all}\,\,x\in G,\,\,t\in[0,T],$$
 we can write
 $$
\int_0^\tau\sum_{x\in G}\left[\rho(x)\eta(x)\partial_t \theta(x,t)+\Delta (\eta(x) \theta(x,t))\right]|u(x,t)|^p \mu(x)\,dt\ge\sum_{x\in G}\rho(x)|u(x,\tau)|^p\eta(x) \theta(x,\tau)\,\mu(x).
 $$
Observe that our choice of $v$ has finite support due to the presence of $\eta$. By using Lemma \ref{lemma5} we compute
 $$
 \Delta (\eta(x) \theta(x,t))=\eta(x)\Delta \theta(x,t)+ \theta(x,t)\Delta\eta(x)+\frac1{\mu(x)}\sum_{y\in G}\nabla_{x,y} \theta\nabla_{xy}\eta\omega(x,y)\,.
 $$
 Therefore we get, by applying also Lemma \ref{lemma51}
 \begin{equation}\label{eq510}
 \begin{aligned}
 \sum_{x\in G}\rho(x)&|u(x,\tau)|^p\eta(x) \theta(x,\tau)\,\mu(x)\\
 &\le \int_0^\tau\sum_{x\in G}\left[\rho(x)\partial_t \theta(x,t)+\Delta \theta(x,t)\right]|u(x,t)|^p \eta(x)\mu(x)\,dt+\int_0^\tau (J_1+J_2)\,dt\\
 &\le \int_0^\tau (J_1+J_2)\,dt\,,
 \end{aligned}
 \end{equation}
where
\begin{align}
&J_1:=\sum_{x\in G}|u(x,t)|^p \theta(x,t)\,\Delta\eta(x)\,\mu(x)\,,\label{eq511}\\
&J_2:=\sum_{x,y\in G}|u(x,t)|^p\nabla_{xy} \theta\,\nabla_{xy}\eta\,\omega(x,y)\,.\label{eq512}
\end{align}
By exploiting the definitions of $\eta$ and $\theta$ given in \eqref{eq316} and \eqref{eq51} respectively, by \eqref{e50f} and assumption \eqref{e1f}, we get
\begin{equation}\label{eq513}
\begin{aligned}
\left| \int_0^\tau J_{1}\,dt \right| &\le\int_0^\tau\sum_{x\in G}|u(x,t)|^p \theta(x,t)\Delta \eta(x)\mu(x)\,dt\\
&\le \frac{C_0}{\delta R}\int_0^\tau\sum_{x\in G}|u(x,t)|^p \theta(x,t)\mu(x)\chi_{\{(1-\delta)R-2s\le d(x,x_0)\le R\}}\,dt\\
&\le \frac{C_0}{\delta R}\int_0^\tau\sum_{x\in G}|u(x,t)|^pe^{-\gamma t}[d(x,x_0)+k]^{-\alpha}\mu(x)\chi_{\{(1-\delta)R-2s\le d(x,x_0)\le R\}}\,dt\\
&\le \frac{C_0}{\delta R}\left[(1-\delta)R-2s+k\right]^{-\alpha}\int_0^\tau\sum_{x\in B_R(x_0)}|u(x,t)|^p \mu(x)\,dt\,.
\end{aligned}
\end{equation}
Let us now consider $J_2$. Implementing the same reasoning used in Lemma \ref{lemma62}, we obtain
$$
\begin{aligned}
|\nabla_{xy}\theta|&\leq e^{-\gamma t}\alpha d(x,y) [d(x,x_0) + k - s]^{-\alpha -1}\\
  &  \leq e^{-\gamma t}\alpha d(x,y)C_3 [d(x,x_0) + k]^{-\alpha} \hspace{1 cm} \forall x \in G,
\end{aligned}
$$
with $C_3 := \frac{k^\alpha}{(k-s)^{\alpha+1}}.$
Then, by means of Lemma \ref{lemma61} and implementing the latter, \eqref{e14f} and \eqref{e1f} into \eqref{eq512}, we obtain
\begin{equation}\label{eq515}
\begin{aligned}
&\left|\int_0^\tau J_{2}\,dt\right| \le \int_0^\tau \sum_{x, y\in G} |u(x,t)|^p |\nabla_{xy}\theta| |\nabla_{xy}\eta| \omega(x,y)\,dt\\
&\leq \frac{\alpha C_3}{\delta R} \int_0^\tau e^{-\gamma t}\sum_{x,y\in G} |u(x,t)|^p d^2(x,y)[d(x,x_0) + k]^{-\alpha} \chi_{\{(1-\delta)R-2s\le d(x,x_0)\le R\}}\omega(x,y)\,dt \\
&\leq \frac{ \alpha C_3}{\delta R}  \int_0^\tau  \sum_{x\in G} [d(x,x_0) + k]^{-\alpha} |u(x,t)|^p \chi_{\{(1-\delta)R-2s\le d(x,x_0)\le R\}} \sum_{y\in G} d^2(x,y)\omega(x,y)\,dt\\
&\leq \frac{  C_0\alpha C_3}{\delta R}  \int_0^\tau  \sum_{x\in G} [d(x,x_0) + k]^{-\alpha} |u(x,t)|^p \chi_{\{(1-\delta)R-2s\le d(x,x_0)\le R\}} \mu(x)\,dt\\
&\leq \frac{ C_0 \alpha C_3}{\delta R} [(1-\delta)R-2s+k]^{-\alpha} \int_0^\tau  \sum_{x\in B_R(x_0)} |u(x,t)|^p  \mu(x)\,dt.
\end{aligned}
\end{equation}
Let $R_n$ be as in the assumptions of Theorem \ref{teo4}. 
Then by condition \eqref{maincondition3}, we have
$$
\int_0^\tau\sum_{x\in B_{R_n}(x_0)}|u(x,t)|^p \mu(x)\,dt\le C\,[R_n+k]^{\beta}.
$$
Hence, \eqref{eq513} and \eqref{eq515} yield
\begin{equation}\label{rhs510}
\begin{aligned}
&\left| \int_0^\tau J_{1}\,dt \right|\le \frac{C_0}{\delta R_n}\left[(1-\delta)R_n-2s+k\right]^{-\alpha}C\,[R_n+k]^{\beta}\\
&\left| \int_0^\tau J_{2}\,dt \right|\le\frac{ C_0 \alpha C_3}{\delta R_n} [(1-\delta)R_n-2s+k]^{-\alpha} C\,[R_n+k]^{\beta}
\end{aligned}
\end{equation}
Combining together \eqref{eq510} and \eqref{rhs510} we get
\begin{equation}\label{eq516}
 \sum_{x\in G}\rho(x)|u(x,\tau)|^p\eta(x) \theta(x,\tau)\,\mu(x)\le \frac{C_0C}{\delta R_n}(1+\alpha C_3)[(1-\delta)R_n-2s+k]^{-\alpha}[R_n+k]^{\beta}
\end{equation}
We now take the limit as $n\to+\infty$ in \eqref{eq516}. Since $\alpha>\max\{0,\beta-1\}$, $R_n\to+\infty$ as $n\to+\infty$ and $\eta\to1$ as $R_n\to\infty$, the latter gives
$$
 \sum_{x\in G}\rho(x)|u(x,\tau)|^p\theta(x,\tau)\,\mu(x)\le0\,.
$$
Due to \eqref{rhobound3}, since $\theta>0$ in $S_T$ and $|u|^p\ge0$, the latter yield that $u(x,\tau)=0$ for all $x\in G$. Note that $\tau$ is arbitrarily chosen in $(0,T]$ hence the theorem follows.
\end{proof}

\bigskip
\bigskip

\noindent{\bf Acknowledgement} 
The author is funded by the Deutsche Forschungsgemeinschaft (DFG, German Research Foundation) - SFB 1283/2 2021 - 317210226. I thank the anonymous Referee for their comments and valuable suggestions.

\noindent{\bf Conflict of interest}. The author states no conflict of interest.

\noindent{\bf Data availability statement}. There are no data associated with this research.

\end{document}